\documentclass[11pt]{amsart}

\usepackage{amsmath}
\usepackage{amssymb}
\usepackage{amsfonts}
\usepackage{geometry}
\usepackage{url}
\usepackage{color}
\usepackage{blkarray}
\usepackage{comment}
\usepackage{blindtext,rotating}

\usepackage{amsmath}
\usepackage{amssymb}
\usepackage{amsfonts}
\usepackage{geometry}
\usepackage{url}
\usepackage{setspace}
\usepackage{amsthm}
\usepackage[all,cmtip]{xy}
\usepackage{amsmath,amscd}
\usepackage{fourier-orns}

\usepackage{amsthm}
\theoremstyle{definition}
\newtheorem{theorem}{Theorem}[section]
\newtheorem{lemma}[theorem]{Lemma}

\newcommand{\hstar}{\mathfrak{h}^*}
\newcommand{\oh}{\mathcal{O}}
\newcommand{\el}{\mathrm{L}}
\newcommand{\m}{\mathrm{M}}
\newcommand{\mbf}{\mathbf}
\newcommand{\s}{\mathfrak{sl}_2}
\newcommand{\B}{\mathcal{B}}

\newcommand{\Ind}{\operatorname{Ind}}

\newcommand{\Hom}{\operatorname{Hom}}

\newcommand{\Res}{\operatorname{Res}}
\newcommand{\Supp}{\operatorname{Supp}}
\newcommand{\Irr}{\operatorname{Irr}}
\newcommand{\rk}{\operatorname{rk}}
\newcommand{\Ext}{\operatorname{Ext}}
\newcommand{\codim}{\operatorname{codim}}

\begin{document}

\title{Irreducible representations of rational Cherednik algebras for exceptional Coxeter groups, part II: some decomposition matrices of $H_c(E_8)$ and $H_c(F_4)$}

\author{Emily Norton}

\begin{abstract} This paper contains the decomposition matrices for blocks of defect at most $2$ in Category $\mathcal{O}_c(W)$ of the rational Cherednik algebra when $W=E_8$ or $F_4$ with equal parameters $c=1/d$, $d>2$ a regular number of $W$. A corollary of the result is a classification of the dimensions of support of the irreducible modules $\el(\tau)$ in $\mathcal{O}_{1/d}(W)$ except in the following cases:  $W=E_8$, $d=4$ or $6$ and $\el(\tau)$ is in the principal block, or $d=2$ or $3$; $W=F_4$, $d=2$. In particular, this classifies the finite-dimensional modules of $H_c(E_8)$ when $d\neq 2,3,4,6$.
\end{abstract}

\maketitle

\section*{Introduction and methods}

As the title suggests, this is the second installment but perhaps not the last of a project intending to classify irreducible representations in Category $\mathcal{O}$ of rational Cherednik algebras for the exceptional Coxeter groups. ``Classify" in the sense of: calculate the decomposition matrix for Category $\mathcal{O}$, that is, the matrix with rows labeled by Verma modules $\m$ and columns by irreducible modules $\el$, and whose entries are the multiplicities of the irreducibles in the Vermas. The inverse of this matrix encodes the characters of all the irreducible representations in $\mathcal{O}$. The characters identify which irreducible representations are finite-dimensional and what their dimensions are, as well as what the dimension of the support variety of each irreducible representation is.

Complete answers to these questions were given by Chmutova for the dihedral groups in \cite{Ch} and by Balagovic-Puranik for $H_3$ in \cite{BaP}. Part I of this project \cite{N} answered these questions for $H_4$, $E_6$, and $E_7$ at parameters $c=1/d$ except for when $d=2$ -- a complete classification up to equivalences of categories when $c\notin\frac{1}{2}\mathbb{Z}$ -- and calculated the finite-dimensional simple modules in $\mathcal{O}_{1/d}(F_4)$ at equal parameters for $d\neq 2$. In this paper we study decomposition matrices in the case of $E_8$ or $F_4$ with equal parameters, when $d>2$. We find the decomposition matrices of all blocks of defect at most $2$, which gives a complete answer for both Coxeter groups when $d>2$ excepting four difficult blocks of $\mathcal{O}_{1/d}(E_8)$: the principal blocks when $d=4$ or $6$ and the two big blocks when $d=3$. Most parameters for $E_8$ are actually quite easy, and we were able to discover the characters and dimensions of many previously unknown finite-dimensional irreducible representations. 

Since Part I appeared, a joint paper by Griffeth, Gusenbauer, Juteau, and Lanini \cite{GGJL} gave a necessary condition for a simple module in $\mathcal{O}_c(W)$ to be finite-dimensional, and their condition has turned out to be often but not always sufficient. The calculations in our paper do not make use of their result because they were done before that paper appeared. While our paper lay sleeping, dreaming of being completed someday to a classification for all blocks of $\mathcal{O}_c(E_8)$, the authors of the forthcoming paper \cite{LSa} sent a draft of their preprint which contains an enumeration of simple modules of given supports in $\mathcal{O}_c(E_8)$ and $\mathcal{O}_c(F_4)$ and a list of the finite-dimensional irreducibles, which they obtained independently by different methods. So here is our Part II woken up suddenly.

\subsection{Notation, terminology, background, methods} More background and definitions are in Part I of this project \cite{N}. 

Let $W=E_8$ or $F_4$ and let $\mathfrak{h}$ be the reflection representation of $W$ and $\hstar$ its dual. Let $S\subset W$ be the set of reflections and let $c:S\rightarrow\mathbb{C}$ be a conjugation-invariant function. Associated to the data of $W$, $\mathfrak{h}$, and $c$ there is a noncommutative algebra called the rational Cherednik algebra $H_c(W)$ which was defined by Etingof and Ginzburg in \cite{EG}. As a vector space, $H_c(W)\cong S\mathfrak{h}\otimes \mathbb{C}[W]\otimes S\hstar$, and $S\mathfrak{h}$, $S\hstar$ are commutative subalgebras of $H_c(W)$. In this paper we only consider $c$ taking the same value on all reflections (this is called ``equal parameters" when there is more than one root-length), i.e. $c\in\mathbb{C}$. We will always take the parameter $c=1/d$ where $d$ is a divisor of a fundamental degree of $W$ and $d>2$ for the following reasons: the category $\mathcal{O}_c(W)$ (defined below) is semi-simple if $c\notin \mathbb{Q}$ with $c=r/d$, $d$ a divisor of a fundamental degree and $\mathrm{gcd}(r,d)=1$; equivalences of categories reduce to the case $c=1/d$; and certain of our methods do not apply when $d=2$.

The Greek letters $\mu,\lambda,\tau,\sigma$ are used for irreducible representations of $W=E_8$ or $W=F_4$. $\Irr W$ denotes the set of irreducible representations of $W$ over $\mathbb{C}$.  $\mathcal{O}_c(W)$ is the category $\mathcal{O}$ of $H_c(W)$ at parameter $c$ defined and studied in \cite{GGOR}. Its objects are all $H_c(W)$-modules which are finitely generated over $S\hstar$ and on which $\mathfrak{h}$ acts locally nilpotently. Basic objects of $\mathcal{O}_c(W)$ are $\el(\tau)$ and $\m(\tau)$ as $\tau$ ranges over $\Irr W$, where $\el(\tau)$ denotes the simple module and $\m(\tau)$ denotes the Verma (or ``standard") module of lowest weight $\tau$. $\mathcal{O}_c(W)$ is a lowest weight category, there is a surjection $\m(\tau)\twoheadrightarrow \el(\tau)$, and $\el(\tau)$ has no self-extensions. A basic problem is to find the multiplicities $[\m(\tau):\el(\tau)]$, the number of times $\el(\tau)$ appears in a Jordan-H\"{o}lder series of $\m(\tau)$. We call $[\m(\tau):\el(\tau)]$ a decomposition number. Since our convention is to use lowest not highest weights, the decomposition number $[\m(\tau):\el(\sigma)]$ can only be nonzero when $\sigma\geq\tau$. Note that $[\m(\tau):\el(\tau)]=1$ always.

By a slight abuse of notation, $H_q(W)$ denotes the Hecke algebra at the parameter $q$, and we will always take $q=e^{2\pi i c}$, where $c$ is the parameter for the Cherednik algebra. The $KZ$ functor $KZ:\mathcal{O}_c(W)\longrightarrow H_q(W)-\mathrm{mod}_{f.d.}$ is a quotient functor and kills all simple modules of less than full support (i.e. those whose support is not all of $\mathfrak{h}$). Here, the support of a simple module $\el(\tau)$ is defined to be $\mathrm{Spec}$ of the annihilator in $S\hstar$ of $\el(\tau)$. The simple modules whose support is the origin are the finite-dimensional modules.

The grading element $h:=\sum_{i=1}^{\rk W} x_iy_i+y_ix_i$ acts by a scalar $h_c(\tau)$ on the lowest weight $\tau$ of Verma or simple module in $\mathcal{O}$. Since $[\m(\tau):\el(\sigma)]>0$ implies that $h_c(\sigma)-h_c(\tau)\in\mathbb{Z}_{\geq0}$ \cite{GGOR}, it is useful to graph the values of $h_c(\tau)$ on a number line as in \cite{BaP}. Such graphs of the $h_c$-weights appear in most of our calculations alongside the decomposition matrices. The element $h_c$ commutes with the action of $W$ and puts a $\mathbb{Z}$-grading on modules in $\mathcal{O}_c(W)$. The graded piece of a module $M$ in degree $n$ will be written as $M[n]$.

Simples $\el(\tau)$ and $\el(\sigma)$ are said to be in the same block of $\mathcal{O}_c(W)$ if there exist $\mu_1=\tau,$ $\mu_2,...,\mu_{k-1},\mu_k=\sigma$ with $\mu_i\in\Irr W$, such that $\Ext^1_{\mathcal{O}_c}(\el(\mu_i),\el(\mu_{i+1}))> 0$ for all $i=1,...,k-1$. A block of $\mathcal{O}_c(W)$ may be considered by shorthand as a subset of $\Irr W$, that is, we may write that $\tau$ is in the block instead of writing $\el(\tau)$ is in the block. The blocks of $\mathcal{O}_c(W)$ then partition $\Irr W$. By ``the principal block" we mean the block containing $\mathrm{Triv}$. Sometimes we will write $\mathcal{B}(\tau)$ for the block of $\mathcal{O}_c(W)$ containing $\tau$. The blocks for $\mathcal{O}_c(W)$ and $H_q(W)-\mathrm{mod}_{f.d.}$ coincide by the Double Centralizer Theorem of \cite{GGOR}. Therefore the first step in finding decomposition numbers for $W$ an exceptional Coxeter group is to look up the blocks for $H_q(W)-\mathrm{mod}_{f.d.}$ in Geck-Pfeiffer's book \cite{GP}.

The notion of the defect of a block comes from group theory; the defect is a number that measures the complexity of the block. It carries over to the Hecke algebras \cite{GP} and thus by the Double Centralizer Theorem to Cherednik category $\mathcal{O}$. By the defect of a block of the Cherednik algebra we mean the defect of the block with the same characters in $H_q(W)-\mathrm{mod}_{f.d.}.$ The blocks of defect $0$ are those $\tau$ in a block by themselves. The blocks of defect $1$ are determined by Brauer trees \cite{GJ} and have a quiver description familiar to many Lie theorists as they also appear as the quiver for parabolic category $\mathcal{O}(\mathfrak{g})^{\mathfrak{p}}$ in type $A_{n-1}\subset A_n$. Defect $1$ blocks have only a single irreducible representation of less than full support. It is basically trivial to classify irreducible representations $\el(\tau)$ and their characters and dimensions of support if they fall into a defect $1$ block, because these blocks are listed in \cite{GJ} and the list of characters basically tells everything about the block. We list the supports of the representations in such blocks and encode their characters in Section 1 anyway, in case it is useful to anybody. Defect $2$ blocks are more complicated but usually still tractable. All the blocks studied in this paper (after the list of defect $1$ blocks) are defect $2$ blocks. Higher defect blocks offer greater challenges. In $\mathcal{O}_c(F_4)$ at equal parameters and $d>2$, all nontrivial blocks have defect $1$ or $2$ \cite{GP}. In $\mathcal{O}_c(E_8)$ and $d>2$, all nontrivial blocks have defect $1$ or $2$ except for the principal blocks when $d=3,4,6$ and the block containing the standard representation of $E_8$ when $d=3$ \cite{GP}.

The rows and columns of the decomposition matrix of a block are labeled by the characters of $\Irr W$ belonging to the block, in nondecreasing order with respect to the values $h_c$ takes on them. The characters $\tau$ with a $\star$ to their left are those such that $\el(\tau)$ is finite-dimensional. Then $\dim\Supp\el(\tau)=0$ and the support is a single point, hence the notation of $\star$. If $\dim\Supp\el(\tau)=k>0$ and $k<\rk W$ then $(k)$ appears to the left of $\tau$ by the decomposition matrix. The rows of the decomposition matrix with no $(k)$ or $\star$ to their left are labeled by $\tau$ such that $\el(\tau)$ has full support. The corresponding columns labeled $\tau$ are the columns of the decomposition matrix of $H_q(W)$. This is the part of the decomposition matrix that is gotten for free, from \cite{GJ}: the decomposition matrix of the Hecke algebra $H_q(W)$ is a submatrix of the decomposition matrix of $\mathcal{O}_c(W)$. Therefore the second step after looking up the blocks in \cite{GP} is to import the decomposition matrix of a block of $H_q(W)$ from Geck-Jacon's book \cite{GJ} as the starting point for the decomposition matrix of the corresponding block of $\mathcal{O}_c(W)$. (We mean the same thing by ``the decomposition matrix of $H_c(W)$" and ``the decomposition matrix of $\mathcal{O}_c(W)$.")

The third step of filling out the decomposition matrix is to apply Lemmas 3.5 and 3.6 from Part I of this project \cite{N}:
\begin{lemma}\label{dim Hom} Suppose $c$ is not a half-integer. Then $$\dim\Hom(\m(\sigma),\m(\tau))=\dim\Hom(\m(\tau'),\m(\sigma'))$$
where $\lambda'=\lambda\otimes\mathrm{Sign}$ for $\lambda\in\Irr W$.
\end{lemma}
\begin{lemma}\label{RR} Let $\tau<_c\sigma$. Suppose there does not exist a chain $\tau=:\nu_{k+1}<_c\nu_k<_c...<_c\nu_2<_c\nu_1<_c\nu_0:=\sigma$ for $k>0$ with $[\m(\nu_i),\el(\nu_{i-1})]>0 $ for all $i=1,\dots,k+1$. Then $[\m(\tau):\el(\sigma)]=\dim\Hom(\m(\sigma),\m(\tau)).$
\end{lemma}
For a number of the defect $2$ decomposition matrices of $\mathcal{O}_c(E_8)$, these lemmas applied to the rectangular decomposition matrices of $H_q(W)$ from \cite{GJ} are enough to finish the job.

For decomposition matrices that can't be completed by the two lemmas, it's necessary to play with induction and restriction from categories $\mathcal{O}$ for parabolic subgroups $W'\subset W$: as defined by Bezrukavnikov-Etingof \cite{BE}, there are exact \cite{BE} and biadjoint (\cite{BE},\cite{Sh},\cite{Lo}) functors $\Res^{\mathcal{O}_c(W)}_{\mathcal{O}_c(W')}:\mathcal{O}_c(W)\longrightarrow\mathcal{O}_c(W')$ and $\Ind^{\mathcal{O}_c(W)}_{\mathcal{O}_c(W')}:\mathcal{O}_c(W')\longrightarrow\Res^{\mathcal{O}_c(W)}_{\mathcal{O}_c(W')}:\mathcal{O}_c(W)$. Sometimes we write $\Res$ for $\Res^{\mathcal{O}_c(W)}_{\mathcal{O}_c(W')}$ when the context is clear. Likewise with $\Ind$. In the Grothendieck group of $\mathcal{O}_c(W)$, $\Res$ and $\Ind$ act the same on a Verma module $\m(\tau)$ as the group representation theoretic induction and restriction do on $\tau$; that is, if $Res^W_{W'}\tau=\sigma_1\oplus\dots\oplus\sigma_k$ then $[\Res^{\mathcal{O}_c(W)}_{\mathcal{O}_c(W')}\m(\tau)]=[\m(\sigma_1)]+\dots+[\m(\sigma_k)]$ \cite{BE}. We may write $\tau\uparrow$ for $\Ind_{W'}^W\tau$ and $\tau\downarrow$ for $\Res_{W'}^W\tau$ when the context is clear. We write $\Res|_{\mathcal{B}}$ for the functor of restriction followed by projection to the block $\mathcal{B}$. We abuse notation and write $\Res$ and $\Ind$ both for the restriction functors on the categories and on the Grothendieck groups. We then write things like $\Res \el=\m_1-\m2+...$ as a shorthand for $[\Res \el]=[\m_1]-[\m_2]+...$ in the Grothendieck group.

The graded character of $\el(\tau)$ is found from row $\tau$ of the inverse matrix $A=(a_{\lambda,\mu})$ to the decomposition matrix, via the formula $\chi_{\el(\tau)}(t)=\left(\sum_{\sigma}a_{\tau,\sigma}\dim_{\mathbb{C}}(\sigma) t^{h_c(\sigma}\right)/(1-t)^{\rk W}$. The order of its pole at $t=1$ equals $\dim\Supp\el(\tau)$. When $\chi_{\el(\tau)}(t)$ is a polynomial then $\el(\tau)$ is finite-dimensional with dimension equal to $\chi_{\el(\tau)}(1)$.

\subsection{Acknowledgments} Thanks to Olivier Dudas and Pavel Etingof for helpful conversations during the time I was computing these decomposition matrices. Thanks to Seth Shelley-Abrahamson and Ivan Losev for communicating their results to me before making them public. 

\section{Defect 1 blocks of $\mathcal{O}_c(E_8)$}

Defect $1$ blocks were explained by Rouquier in \cite{R}, Theorem 5.15. The decomposition matrix of a defect $1$ block always has $1$'s on and just above the diagonal, and $0$'s everywhere else. Thus if $\mathcal{A}$ is a defect $1$ block, it suffices to know the characters of $\Irr W$ labeling irreducibles in $\mathcal{A}$ together with the order on them determined by increasing values of $h_c$ in order to know the decomposition matrix of $\mathcal{A}$. The information about the characters and their order can be found in the back of \cite{GP}. To write down the characters of the Cherednik algebra modules $L(\lambda)$ for $\lambda\in\mathcal{A}$ it's necessary to actually calculate the values $h_c(\lambda)$. We print this information below for the cases that have not been written down before.

The finite-dimensional simple representations belonging to a defect $1$ block were classified by Rouquier \cite{R}. He calculated their characters for $c=1/d$, which also implies their dimensions (by evaluating the character at $t=1$). This information can be found in \cite{R} following Theorem 5.15. Note that Rouquier's convention for the $h_c$-function differs from ours by a shift -- he normalizes so that on the trivial representation it always takes value $0$.

Remark 5.16 of \cite{R} should generalize to the statement that there is exactly one simple module, say $\el(\lambda_0)$, of less than full support in any defect $1$ block $\mathcal{A}$ of $\mathcal{O}_c(W)$; that $h_c(\lambda_0)$ is minimal over $h_c(\mu)$ for $\mu\in\mathcal{A}$; and that \\$\#\{\el(\mu)\hbox{ simple }|\;\el(\mu)\in\mathcal{A}\}\geq \codim\Supp\el(\lambda_0)+1$. In all examples in this paper, the inequality is an equality. If this is true for a conceptual reason, then to see the dimensions of support of all modules in a defect $1$ block, it would not even be necessary to calculate anything. It would be enough to look at Tables F.1-F.6 in \cite{GP}, and if there are $m$ characters in the block then the dimension of support of $L(\tau_1)$ where $\tau_1$ is the first representation listed will be equal to $\rk W +1-m$.
 
The $h_c$-weights of each defect $1$ block of $\mathcal{O}_c(E_8)$ at $c=1/d$ which does not contain a finite-dimensional representation, and the dimension of support of the unique irreducible representation belonging to each such block, appears below.  All data about the defects of blocks and the characters in the blocks is from Table F.6 of \cite{GP}. For the first few values of $d$ the $h_c$-weights of the $\tau$ in the block are drawn on a number line and labeled by the $\tau$. For the rest of the values of $d$ the $h_c$-weights have been shifted so that the smallest weight in the block is $0$. They are listed in the form $\{h_c'(\tau)\}=\{h_c'(\tau_1),h_c'(\tau_2)...,h_c'(\tau_m)\}$, where $h_c'=h_c-h_c(\tau_1)$. The reader can find in Table F.6 of \cite{GP} the list of the $\tau$ belonging to the block with the specified irreducible $L(\tau_1)$ of less than full support. The dimension of support of $L(\tau_1)$ is the power of $(1-t)$ in the denominator of $\sum_{k=1}^m (-1)^{k-1} \dim\tau_k t^{h_c'(\tau_k)}/(1-t)^8$ (after canceling any factors of $(1-t)$ in the numerator).

\subsection{{\Large$\mbf{c=\frac{1}{18}}$}} 

\small
\begin{center}
\begin{picture}(300,50)
\put(0,30){\line(1,0){300}}
\put(0,30){\circle*{5}}
\put(120,30){\circle*{5}}
\put(150,30){\circle*{5}}
\put(180,30){\circle*{5}}
\put(210,30){\circle*{5}}
\put(240,30){\circle*{5}}
\put(270,30){\circle*{5}}
\put(300,30){\circle*{5}}
\put(-12,40){$-2\frac{2}{3}$}
\put(115,40){$1\frac{1}{3}$}
\put(145,40){$2\frac{1}{3}$}
\put(175,40){$3\frac{1}{3}$}
\put(205,40){$4\frac{1}{3}$}
\put(235,40){$5\frac{1}{3}$}
\put(265,40){$6\frac{1}{3}$}
\put(295,40){$7\frac{1}{3}$}
\put(-4,15){$1_x$}
\put(112,15){$210_x$}
\put(140,15){$1008_z$}
\put(170,15){$2100_x$}
\put(200,15){$2400_z'$}
\put(230,15){$1575_x'$}
\put(262,15){$560_z'$}
\put(295,15){$84_x'$}
\end{picture}

\begin{picture}(300,50)
\put(0,30){\line(1,0){300}}
\put(0,30){\circle*{5}}
\put(30,30){\circle*{5}}
\put(60,30){\circle*{5}}
\put(90,30){\circle*{5}}
\put(120,30){\circle*{5}}
\put(150,30){\circle*{5}}
\put(180,30){\circle*{5}}
\put(300,30){\circle*{5}}
\put(-5,40){$\frac{2}{3}$}
\put(25,40){$1\frac{2}{3}$}
\put(55,40){$2\frac{2}{3}$}
\put(85,40){$3\frac{2}{3}$}
\put(115,40){$4\frac{2}{3}$}
\put(145,40){$5\frac{2}{3}$}
\put(175,40){$6\frac{2}{3}$}
\put(295,40){$10\frac{1}{3}$}
\put(-8,15){$84_x$}
\put(22,15){$560_z$}
\put(50,15){$1575_x$}
\put(80,15){$2400_z$}
\put(110,15){$2100_x'$}
\put(140,15){$1008_z'$}
\put(172,15){$210_x'$}
\put(295,15){$1_x'$}
\end{picture}

\end{center}

\normalsize

\begin{itemize}
\item $\dim\Supp\el(1_x)=1$
\item $\dim\Supp\el(84_x)=1$
\end{itemize}

\subsection{{\Large$\mbf{c=\frac{1}{14}}$}} 

\small
\begin{center}
\begin{picture}(450,50)
\put(0,30){\line(1,0){450}}
\put(0,30){\circle*{5}}
\put(180,30){\circle*{5}}
\put(210,30){\circle*{5}}
\put(240,30){\circle*{5}}
\put(270,30){\circle*{5}}
\put(300,30){\circle*{5}}
\put(360,30){\circle*{5}}
\put(450,30){\circle*{5}}
\put(-8,40){$-4\frac{4}{7}$}
\put(175,40){$1\frac{3}{7}$}
\put(205,40){$2\frac{3}{7}$}
\put(235,40){$3\frac{3}{7}$}
\put(265,40){$4\frac{3}{7}$}
\put(295,40){$5\frac{3}{7}$}
\put(355,40){$7\frac{3}{7}$}
\put(445,40){$10\frac{3}{7}$}
\put(-2,15){$1_x$}
\put(172,15){$700_x$}
\put(200,15){$3240_x$}
\put(230,15){$6075_x$}
\put(260,15){$5600_z'$}
\put(290,15){$2268_x'$}
\put(352,15){$210_x'$}
\put(445,15){$8_z'$}
\end{picture}

\begin{picture}(450,50)
\put(0,30){\line(1,0){450}}
\put(0,30){\circle*{5}}
\put(90,30){\circle*{5}}
\put(150,30){\circle*{5}}
\put(180,30){\circle*{5}}
\put(210,30){\circle*{5}}
\put(240,30){\circle*{5}}
\put(270,30){\circle*{5}}
\put(450,30){\circle*{5}}
\put(-8,40){$-2\frac{3}{7}$}
\put(85,40){$\frac{4}{7}$}
\put(145,40){$2\frac{4}{7}$}
\put(175,40){$3\frac{4}{7}$}
\put(205,40){$4\frac{4}{7}$}
\put(235,40){$5\frac{4}{7}$}
\put(265,40){$6\frac{4}{7}$}
\put(445,40){$10\frac{4}{7}$}
\put(-2,15){$8_z$}
\put(82,15){$210_x$}
\put(140,15){$2268_x$}
\put(170,15){$5600_z$}
\put(200,15){$6075_x'$}
\put(230,15){$3240_x'$}
\put(262,15){$700_x'$}
\put(445,15){$1_x'$}
\end{picture}

\end{center}
\normalsize

\begin{itemize}
\item$\dim\Supp\el(1_x)=1$
\item$\dim\Supp\el(8_z)=1$
\end{itemize}

\subsection{{\Large$\mbf{c=\frac{1}{12}}$}}
\small
\begin{center}
\begin{picture}(300,50)
\put(0,30){\line(1,0){300}}
\put(0,30){\circle*{5}}
\put(120,30){\circle*{5}}
\put(150,30){\circle*{5}}
\put(210,30){\circle*{5}}
\put(240,30){\circle*{5}}
\put(270,30){\circle*{5}}
\put(300,30){\circle*{5}}
\put(-12,40){$-3.5$}
\put(115,40){$.5$}
\put(145,40){$1.5$}
\put(205,40){$3.5$}
\put(235,40){$4.5$}
\put(265,40){$5;5$}
\put(295,40){$6.5$}
\put(-4,15){$8_z$}
\put(112,15){$560_z$}
\put(140,15){$1344_x$}
\put(200,15){$3200_x$}
\put(230,15){$4200_x'$}
\put(262,15){$2240_x'$}
\put(295,15){$448_z'$}
\end{picture}

\begin{picture}(300,50)
\put(0,30){\line(1,0){300}}
\put(0,30){\circle*{5}}
\put(30,30){\circle*{5}}
\put(60,30){\circle*{5}}
\put(90,30){\circle*{5}}
\put(150,30){\circle*{5}}
\put(180,30){\circle*{5}}
\put(300,30){\circle*{5}}
\put(-5,40){$1.5$}
\put(25,40){$2.5$}
\put(55,40){$3.5$}
\put(85,40){$4.5$}
\put(145,40){$6.5$}
\put(175,40){$7.5$}
\put(295,40){$11.5$}
\put(-8,15){$448_z$}
\put(22,15){$2240_x$}
\put(50,15){$4200_x$}
\put(80,15){$3200_x'$}
\put(140,15){$1344_x'$}
\put(172,15){$560_z'$}
\put(295,15){$8_z'$}
\end{picture}

\begin{picture}(300,50)
\put(0,30){\line(1,0){300}}
\put(0,30){\circle*{5}}
\put(60,30){\circle*{5}}
\put(120,30){\circle*{5}}
\put(150,30){\circle*{5}}
\put(180,30){\circle*{5}}
\put(240,30){\circle*{5}}
\put(300,30){\circle*{5}}
\put(-8,40){$-1$}
\put(55,40){$1$}
\put(115,40){$3$}
\put(145,40){$4$}
\put(175,40){$5$}
\put(235,40){$7$}
\put(295,40){$9$}
\put(-6,15){$84_x$}
\put(52,15){$700_x$}
\put(110,15){$4200_x$}
\put(140,15){$7168_w$}
\put(170,15){$4200_x'$}
\put(232,15){$700_x'$}
\put(295,15){$84_x'$}
\end{picture}
\end{center}
\normalsize

\begin{itemize}
\item$\dim\Supp\el(8_z)=2$
\item$\dim\Supp\el(448_z)=2$
\item$\dim\Supp\el(84_x)=2$
\end{itemize}

\subsection{{\Large$\mbf{c=\frac{1}{10}}$}}
\small
\begin{center}
\begin{picture}(360,50)
\put(0,30){\line(1,0){360}}
\put(0,30){\circle*{5}}
\put(90,30){\circle*{5}}
\put(120,30){\circle*{5}}
\put(180,30){\circle*{5}}
\put(210,30){\circle*{5}}
\put(270,30){\circle*{5}}
\put(360,30){\circle*{5}}
\put(-10,40){$-.8$}
\put(84,40){$2.2$}
\put(114,40){$3.2$}
\put(174,40){$5.2$}
\put(204,40){$6.2$}
\put(264,40){$8.2$}
\put(353,40){$11.2$}
\put(-6,15){$50_x$}
\put(78,15){$1400_{zz}$}
\put(110,15){$2835_x$}
\put(170,15){$4200_x'$}
\put(200,15){$3240_z'$}
\put(262,15){$560_z'$}
\put(357,15){$35_x'$}
\end{picture}

\begin{picture}(360,50)
\put(0,30){\line(1,0){360}}
\put(0,30){\circle*{5}}
\put(90,30){\circle*{5}}
\put(150,30){\circle*{5}}
\put(180,30){\circle*{5}}
\put(240,30){\circle*{5}}
\put(270,30){\circle*{5}}
\put(360,30){\circle*{5}}
\put(-12,40){$-3.8$}
\put(80,40){$-.8$}
\put(144,40){$1.2$}
\put(174,40){$2.2$}
\put(234,40){$4.2$}
\put(264,40){$5.2$}
\put(353,40){$8.2$}
\put(-6,15){$35_x$}
\put(78,15){$560_z$}
\put(140,15){$3240_z$}
\put(170,15){$2835_x$}
\put(230,15){$4200_x'$}
\put(262,15){$1400_{zz}'$}
\put(357,15){$50_x'$}
\end{picture}
\end{center}
\normalsize

\begin{itemize}
\item $\dim\Supp\el(50_x)=2$
\item$\dim\Supp\el(35_x)=2$
\end{itemize}

\subsection{{\Large$\mbf{c=\frac{1}{9}}$}}

\begin{itemize}
\item $\dim\Supp\el(1_x)=2$\qquad\qquad$\{h_c'(\tau)\}=\{0,8,12,14,15,16,20\}$
\item $\dim\Supp\el(28_x)=2$\qquad\qquad$\{h_c'(\tau)\}=\{0,4,5,6,8,12,20\}$
\item $\dim\Supp\el(160_z)=2$\qquad\qquad$\{h_c'(\tau)\}=\{0,2,4,6,7,10,12\}$
\item $\dim\Supp\el(112_z)=2$\qquad\qquad$\{h_c'(\tau)\}=\{0,2,3,4,6,8,10\}$
\item $\dim\Supp\el(8_z)=2$\qquad\qquad$\{h_c'(\tau)\}=\{0,6,8,10,12,14,20\}$
\item $\dim\Supp\el(35_x)=2$\qquad\qquad$\{h_c'(\tau)\}=\{0,4,6,8,10,12,16\}$
\end{itemize}

\subsection{{\Large$\mbf{c=\frac{1}{8}}$}}

\begin{itemize}
\item $\dim\Supp\el(8_z)=3$\qquad\qquad$\{h_c'(\tau)\}=\{0,9,10,12,14,15\}$
\item $\dim\Supp\el(1344_x)=3$\qquad\qquad$\{h_c'(\tau)\}=\{0,1,3,5,6,15\}$
\item $\dim\Supp\el(28_x)=3$\qquad\qquad$\{h_c'(\tau)\}=\{0,5,6,9,10,12\}$
\item $\dim\Supp\el(300_x)=3$\qquad\qquad$\{h_c'(\tau)\}=\{0,2,3,6,7,12\}$
\item $\dim\Supp\el(112_z)=3$\qquad\qquad$\{h_c'(\tau)\}=\{0,3,5,9,12,15\}$
\item $\dim\Supp\el(84_x)=3$\qquad\qquad$\{h_c'(\tau)\}=\{0,3,6,10,12,15\}$
\end{itemize}

\subsection{{\Large$\mbf{c=\frac{1}{7}}$}}
\begin{itemize}
\item $\dim\Supp\el(1_x)=2$\qquad\qquad $\{h_c'(\tau)\}=\{0,12,15,16,18,20,24\}$
\item $\dim\Supp\el(50_x)=2$\qquad\qquad$\{h_c'(\tau)\}=\{0,4,6,8,9,12,24\}$
\item $\dim\Supp\el(8_z)=2$\qquad\qquad $\{h_c'(\tau)\}=\{0,6,10,12,15,16,18\}$
\item $\dim\Supp\el(400_z)=2$\qquad\qquad$\{h_c'(\tau)\}=\{0,2,3,6,8,12,18\}$
\end{itemize}

\subsection{{\Large$\mbf{c=\frac{1}{6}}$}}
\begin{itemize}
\item $\dim\Supp\el(972_x)=4$\qquad\qquad $\{h_c'(\tau)\}=\{0,2,5,7,10\}$
\item$\dim\Supp\el(567_x)=4$\qquad\qquad $\{h_c'(\tau)\}=\{0,3,5,8,10\}$
\item$\dim\Supp\el(1008_z)=4$\qquad\qquad $\{h_c'(\tau)\}=\{0,1,5,9,10\}$
\end{itemize}

\subsection{{\Large$\mbf{c=\frac{1}{5}}$}}
\begin{itemize}
\item $\dim\Supp\el(210_x)=4$\qquad\qquad $\{h_c'(\tau)\}=\{0,6,8,12,24\}$
\item $\dim\Supp\el(35_x)=4$\qquad\qquad $\{h_c'(\tau)\}=\{0,12,16,18,24\}$
\item $\dim\Supp\el(160_z)=4$\qquad\qquad $\{h_c'(\tau)\}=\{0,6,12,14,18\}$
\item $\dim\Supp\el(560_z)=4$\qquad\qquad $\{h_c'(\tau)\}=\{0,4,6,12,18\}$
\end{itemize}

\subsection{{\Large$\mbf{c=\frac{1}{3}}$}}
\begin{itemize}
\item $\dim\Supp\el(567_x)=6$\qquad\qquad$\{h_c'(\tau)\}=\{0,16,20\}$
\item $\dim\Supp\el(2268_x)=6$\qquad\qquad$\{h_c'(\tau)\}=\{0,4,20\}$
\item $\dim\Supp\el(1296_z)=6$\qquad\qquad$\{h_c'(\tau)\}=\{0,10,14\}$
\item $\dim\Supp\el(3240_z)=6$\qquad\qquad$\{h_c'(\tau)\}=\{0,4,14\}$
\item $\dim\Supp\el(1008_z)=6$\qquad\qquad$\{h_c'(\tau)\}=\{0,10,20\}$
\item $\dim\Supp\el(1575_x)=6$\qquad\qquad$\{h_c'(\tau)\}=\{0,8,16\}$
\end{itemize}

\section{Defect $2$ blocks of $\mathcal{O}_c(E_8)$}
The decomposition matrices for $c=\frac{1}{12}$, $\frac{1}{10},$ $\frac{1}{8}$, $\frac{1}{6}$, $\frac{1}{5}$, and $\frac{1}{4}$ which appear in this section are almost entirely given by the decomposition matrices of the Hecke algebra $H_q(E_8)$ at $q=e^{2\pi i c}$, which are in \cite{GJ}.  In the case of $c=\frac{1}{12},\frac{1}{10},$ $\frac{1}{8}$, $\frac{1}{5}$, there are six additional columns which have to be found: but these follow immediately from Lemmas \ref{RR} and \ref{dim Hom}. Thus the decomposition matrices at these parameters are printed without any proofs below. In the case of $c=\frac{1}{6}$ and $\frac{1}{4}$, there are five additional columns to find, and in these cases one has to make additional arguments to complete these columns.

\subsection{{\Large$\mbf{c=\frac{1}{12}}$}}

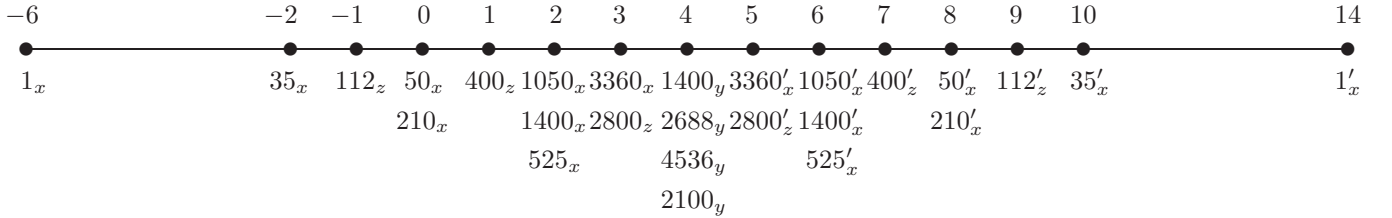
\small
\begin{center}
\begin{picture}(500,100)
\put(0,80){\line(1,0){500}}
\put(0,80){\circle*{5}}
\put(100,80){\circle*{5}}
\put(125,80){\circle*{5}}
\put(150,80){\circle*{5}}
\put(175,80){\circle*{5}}
\put(200,80){\circle*{5}}
\put(225,80){\circle*{5}}
\put(250,80){\circle*{5}}
\put(275,80){\circle*{5}}
\put(300,80){\circle*{5}}
\put(325,80){\circle*{5}}
\put(350,80){\circle*{5}}
\put(375,80){\circle*{5}}
\put(400,80){\circle*{5}}
\put(500,80){\circle*{5}}
\put(-8,90){$-6$}
\put(90,90){$-2$}
\put(115,90){$-1$}
\put(148,90){$0$}
\put(173,90){$1$}
\put(198,90){$2$}
\put(222,90){$3$}
\put(247,90){$4$}
\put(272,90){$5$}
\put(297,90){$6$}
\put(322,90){$7$}
\put(347,90){$8$}
\put(372,90){$9$}
\put(395,90){$10$}
\put(495,90){$14$}
\put(-2,65){$1_x$}
\put(92,65){$35_x$}
\put(117,65){$112_z$}
\put(143,65){$50_x$}
\put(140,50){$210_x$}
\put(166,65){$400_z$}
\put(187,65){$1050_x$}
\put(187,50){$1400_x$}
\put(190,35){$525_x$}
\put(213,65){$3360_x$}
\put(213,50){$2800_z$}
\put(240,65){$1400_y$}
\put(240,50){$2688_y$}
\put(240,35){$4536_y$}
\put(240,20){$2100_y$}
\put(266,65){$3360_x'$}
\put(266,50){$2800_z'$}
\put(292,65){$1050_x'$}
\put(292,50){$1400_x'$}
\put(295,35){$525_x'$}
\put(318,65){$400_z'$}
\put(345,65){$50_x'$}
\put(342,50){$210_x'$}
\put(367,65){$112_z'$}
\put(395,65){$35_x'$}
\put(495,65){$1_x'$}
\end{picture}
\end{center}

\small
\[
\begin{blockarray}{cccccccccccccccccccccccccccc}
\begin{block}{cc(cccccccccccccccccccccccccc)}
\star&1_x&1&\cdot&1&\cdot&\cdot&1&\cdot&1&\cdot&\cdot&\cdot&\cdot&\cdot&\cdot&\cdot&\cdot&\cdot&\cdot&\cdot&\cdot&\cdot&\cdot&\cdot&\cdot&\cdot&\cdot\\
\star&35_x&\cdot&1&1&1&\cdot&\cdot&\cdot&\cdot&\cdot&\cdot&\cdot&\cdot&\cdot&\cdot&\cdot&\cdot&\cdot&\cdot&\cdot&\cdot&\cdot&\cdot&\cdot&\cdot&\cdot&\cdot\\
\texttt{(2)}&112_z&\cdot&\cdot&1&1&\cdot&\cdot&1&1&\cdot&1&\cdot&\cdot&\cdot&\cdot&\cdot&\cdot&\cdot&\cdot&\cdot&\cdot&\cdot&\cdot&\cdot&\cdot&\cdot&\cdot\\
&210_x&\cdot&\cdot&\cdot&1&\cdot&\cdot&\cdot&\cdot&\cdot&1&\cdot&\cdot&\cdot&\cdot&1&\cdot&\cdot&\cdot&\cdot&\cdot&\cdot&\cdot&\cdot&\cdot&\cdot&\cdot\\
\star&50_x&\cdot&\cdot&\cdot&\cdot&1&1&\cdot&\cdot&1&\cdot&\cdot&\cdot&\cdot&\cdot&\cdot&\cdot&\cdot&\cdot&\cdot&\cdot&\cdot&\cdot&\cdot&\cdot&\cdot&\cdot\\
\texttt{(1)}&400_z&\cdot&\cdot&\cdot&\cdot&\cdot&1&\cdot&1&1&\cdot&1&\cdot&\cdot&\cdot&\cdot&\cdot&\cdot&\cdot&\cdot&\cdot&\cdot&\cdot&\cdot&\cdot&\cdot&\cdot\\
\texttt{(2)}&525_x&\cdot&\cdot&\cdot&\cdot&\cdot&\cdot&1&\cdot&\cdot&1&\cdot&1&\cdot&\cdot&\cdot&\cdot&\cdot&\cdot&\cdot&\cdot&\cdot&\cdot&\cdot&\cdot&\cdot&\cdot\\
&1400_x&\cdot&\cdot&\cdot&\cdot&\cdot&\cdot&\cdot&1&\cdot&1&1&\cdot&1&\cdot&\cdot&\cdot&\cdot&\cdot&\cdot&\cdot&\cdot&\cdot&\cdot&\cdot&\cdot&\cdot\\
&1050_x&\cdot&\cdot&\cdot&\cdot&\cdot&\cdot&\cdot&\cdot&1&\cdot&1&\cdot&\cdot&1&\cdot&\cdot&\cdot&\cdot&\cdot&\cdot&\cdot&\cdot&\cdot&\cdot&\cdot&\cdot\\
&2800_z&\cdot&\cdot&\cdot&\cdot&\cdot&\cdot&\cdot&\cdot&\cdot&1&\cdot&1&1&\cdot&1&1&\cdot&\cdot&\cdot&\cdot&\cdot&\cdot&\cdot&\cdot&\cdot&\cdot\\
&3360_z&\cdot&\cdot&\cdot&\cdot&\cdot&\cdot&\cdot&\cdot&\cdot&\cdot&1&\cdot&1&1&\cdot&\cdot&1&\cdot&\cdot&\cdot&\cdot&\cdot&\cdot&\cdot&\cdot&\cdot\\
&2100_y&\cdot&\cdot&\cdot&\cdot&\cdot&\cdot&\cdot&\cdot&\cdot&\cdot&\cdot&1&\cdot&\cdot&\cdot&1&\cdot&1&\cdot&\cdot&\cdot&\cdot&\cdot&\cdot&\cdot&\cdot\\
&4536_y&\cdot&\cdot&\cdot&\cdot&\cdot&\cdot&\cdot&\cdot&\cdot&\cdot&\cdot&\cdot&1&\cdot&\cdot&1&1&\cdot&1&\cdot&\cdot&\cdot&\cdot&\cdot&\cdot&\cdot\\
&2688_y&\cdot&\cdot&\cdot&\cdot&\cdot&\cdot&\cdot&\cdot&\cdot&\cdot&\cdot&\cdot&\cdot&1&\cdot&\cdot&1&\cdot&\cdot&1&\cdot&\cdot&\cdot&\cdot&\cdot&\cdot\\
&1400_y&\cdot&\cdot&\cdot&\cdot&\cdot&\cdot&\cdot&\cdot&\cdot&\cdot&\cdot&\cdot&\cdot&\cdot&1&1&\cdot&\cdot&\cdot&\cdot&\cdot&1&\cdot&\cdot&\cdot&\cdot\\
&2800_z'&\cdot&\cdot&\cdot&\cdot&\cdot&\cdot&\cdot&\cdot&\cdot&\cdot&\cdot&\cdot&\cdot&\cdot&\cdot&1&\cdot&1&1&\cdot&\cdot&1&\cdot&1&\cdot&\cdot\\
&3360_z'&\cdot&\cdot&\cdot&\cdot&\cdot&\cdot&\cdot&\cdot&\cdot&\cdot&\cdot&\cdot&\cdot&\cdot&\cdot&\cdot&1&\cdot&1&1&1&\cdot&\cdot&\cdot&\cdot&\cdot\\
&525_x'&\cdot&\cdot&\cdot&\cdot&\cdot&\cdot&\cdot&\cdot&\cdot&\cdot&\cdot&\cdot&\cdot&\cdot&\cdot&\cdot&\cdot&1&\cdot&\cdot&\cdot&\cdot&\cdot&1&\cdot&\cdot\\
&1400_x'&\cdot&\cdot&\cdot&\cdot&\cdot&\cdot&\cdot&\cdot&\cdot&\cdot&\cdot&\cdot&\cdot&\cdot&\cdot&\cdot&\cdot&\cdot&1&\cdot&1&\cdot&\cdot&1&\cdot&1\\
&1050_x'&\cdot&\cdot&\cdot&\cdot&\cdot&\cdot&\cdot&\cdot&\cdot&\cdot&\cdot&\cdot&\cdot&\cdot&\cdot&\cdot&\cdot&\cdot&\cdot&1&1&\cdot&1&\cdot&\cdot&\cdot\\
&400_z'&\cdot&\cdot&\cdot&\cdot&\cdot&\cdot&\cdot&\cdot&\cdot&\cdot&\cdot&\cdot&\cdot&\cdot&\cdot&\cdot&\cdot&\cdot&\cdot&\cdot&1&\cdot&1&\cdot&\cdot&1\\
&210_x'&\cdot&\cdot&\cdot&\cdot&\cdot&\cdot&\cdot&\cdot&\cdot&\cdot&\cdot&\cdot&\cdot&\cdot&\cdot&\cdot&\cdot&\cdot&\cdot&\cdot&\cdot&1&\cdot&1&1&\cdot\\
&50_x'&\cdot&\cdot&\cdot&\cdot&\cdot&\cdot&\cdot&\cdot&\cdot&\cdot&\cdot&\cdot&\cdot&\cdot&\cdot&\cdot&\cdot&\cdot&\cdot&\cdot&\cdot&\cdot&1&\cdot&\cdot&\cdot\\
&112_z'&\cdot&\cdot&\cdot&\cdot&\cdot&\cdot&\cdot&\cdot&\cdot&\cdot&\cdot&\cdot&\cdot&\cdot&\cdot&\cdot&\cdot&\cdot&\cdot&\cdot&\cdot&\cdot&\cdot&1&1&1\\
&35_x'&\cdot&\cdot&\cdot&\cdot&\cdot&\cdot&\cdot&\cdot&\cdot&\cdot&\cdot&\cdot&\cdot&\cdot&\cdot&\cdot&\cdot&\cdot&\cdot&\cdot&\cdot&\cdot&\cdot&\cdot&1&\cdot\\
&1_x'&\cdot&\cdot&\cdot&\cdot&\cdot&\cdot&\cdot&\cdot&\cdot&\cdot&\cdot&\cdot&\cdot&\cdot&\cdot&\cdot&\cdot&\cdot&\cdot&\cdot&\cdot&\cdot&\cdot&\cdot&\cdot&1\\
\end{block}
\end{blockarray}
\]

\newgeometry{top=2cm, bottom=2cm}
\begin{turn}{270}
\begin{minipage}{\linewidth}
\small
\begin{center}
\[
\begin{blockarray}{cccccccccccccccccccccccccccc}
\begin{block}{cc(cccccccccccccccccccccccccc)}
\star&1_x&1&\cdot&-1&1&\cdot&-1&1&1&1&-2&-1&1&2&\cdot&1&-2&-1&1&1&1&-1&1&\cdot&-1&\cdot&1\\
\star&35_x&\cdot&1&-1&\cdot&\cdot&\cdot&1&1&\cdot&-1&-1&\cdot&1&1&1&-1&-1&1&1&\cdot&\cdot&\cdot&\cdot&-1&1&\cdot\\
\mathtt{(2)}&112_z&\cdot&\cdot&1&-1&\cdot&\cdot&-1&-1&\cdot&2&1&-1&-2&-1&-1&2&2&-1&-2&-1&1&-1&\cdot&2&-1&-1\\
&210_x&\cdot&\cdot&\cdot&1&\cdot&\cdot&\cdot&\cdot&\cdot&-1&\cdot&1&1&\cdot&\cdot&-1&-1&\cdot&1&1&-1&1&\cdot&-1&\cdot&1\\
\star&50_x&\cdot&\cdot&\cdot&\cdot&1&-1&\cdot&1&\cdot&-1&\cdot&1&\cdot&\cdot&1&-1&\cdot&\cdot&1&\cdot&-1&\cdot&1&\cdot&\cdot&\cdot\\
\mathtt{(1)}&400_z&\cdot&\cdot&\cdot&\cdot&\cdot&1&\cdot&-1&-1&1&1&-1&-1&\cdot&-1&2&\cdot&-1&-1&\cdot&1&-1&-1&1&\cdot&-1\\
\mathtt{(2)}&525_x&\cdot&\cdot&\cdot&\cdot&\cdot&\cdot&1&\cdot&\cdot&-1&\cdot&\cdot&1&\cdot&1&-1&-1&1&1&1&-1&\cdot&\cdot&-1&1&1\\
&1400_x&\cdot&\cdot&\cdot&\cdot&\cdot&\cdot&\cdot&1&\cdot&-1&-1&1&1&1&1&-2&-1&1&2&\cdot&-1&1&1&-2&1&1\\
&1050_x&\cdot&\cdot&\cdot&\cdot&\cdot&\cdot&\cdot&\cdot&1&\cdot&-1&\cdot&1&\cdot&\cdot&-1&\cdot&1&\cdot&\cdot&\cdot&1&\cdot&-1&\cdot&1\\
&2800_z&\cdot&\cdot&\cdot&\cdot&\cdot&\cdot&\cdot&\cdot&\cdot&1&\cdot&-1&-1&\cdot&-1&2&1&-1&-2&-1&2&-1&-1&2&-1&-2\\
&3360_z&\cdot&\cdot&\cdot&\cdot&\cdot&\cdot&\cdot&\cdot&\cdot&\cdot&1&\cdot&-1&-1&\cdot&1&1&-1&-1&\cdot&\cdot&-1&\cdot&2&-1&-1\\
&2100_y&\cdot&\cdot&\cdot&\cdot&\cdot&\cdot&\cdot&\cdot&\cdot&\cdot&\cdot&1&\cdot&\cdot&\cdot&-1&\cdot&\cdot&1&\cdot&-1&1&1&-1&\cdot&1\\
&4536_y&\cdot&\cdot&\cdot&\cdot&\cdot&\cdot&\cdot&\cdot&\cdot&\cdot&\cdot&\cdot&1&\cdot&\cdot&-1&-1&1&1&1&-1&1&\cdot&-2&1&2\\
&2688_y&\cdot&\cdot&\cdot&\cdot&\cdot&\cdot&\cdot&\cdot&\cdot&\cdot&\cdot&\cdot&\cdot&1&\cdot&\cdot&-1&\cdot&1&\cdot&\cdot&\cdot&\cdot&-1&1&\cdot\\
&1400_y&\cdot&\cdot&\cdot&\cdot&\cdot&\cdot&\cdot&\cdot&\cdot&\cdot&\cdot&\cdot&\cdot&\cdot&1&-1&\cdot&1&1&\cdot&-1&\cdot&1&-1&1&1\\
&2800_z'&\cdot&\cdot&\cdot&\cdot&\cdot&\cdot&\cdot&\cdot&\cdot&\cdot&\cdot&\cdot&\cdot&\cdot&\cdot&1&\cdot&-1&-1&\cdot&1&-1&-1&2&-1&-2\\
&3360_z'&\cdot&\cdot&\cdot&\cdot&\cdot&\cdot&\cdot&\cdot&\cdot&\cdot&\cdot&\cdot&\cdot&\cdot&\cdot&\cdot&1&\cdot&-1&-1&1&\cdot&\cdot&1&-1&-1\\
&525_x'&\cdot&\cdot&\cdot&\cdot&\cdot&\cdot&\cdot&\cdot&\cdot&\cdot&\cdot&\cdot&\cdot&\cdot&\cdot&\cdot&\cdot&1&\cdot&\cdot&\cdot&\cdot&\cdot&-1&1&1\\
&1400_x'&\cdot&\cdot&\cdot&\cdot&\cdot&\cdot&\cdot&\cdot&\cdot&\cdot&\cdot&\cdot&\cdot&\cdot&\cdot&\cdot&\cdot&\cdot&1&\cdot&-1&\cdot&1&-1&1&1\\
&1050_x'&\cdot&\cdot&\cdot&\cdot&\cdot&\cdot&\cdot&\cdot&\cdot&\cdot&\cdot&\cdot&\cdot&\cdot&\cdot&\cdot&\cdot&\cdot&\cdot&1&-1&\cdot&\cdot&\cdot&\cdot&1\\
&400_z'&\cdot&\cdot&\cdot&\cdot&\cdot&\cdot&\cdot&\cdot&\cdot&\cdot&\cdot&\cdot&\cdot&\cdot&\cdot&\cdot&\cdot&\cdot&\cdot&\cdot&1&\cdot&-1&\cdot&\cdot&-1\\
&210_x'&\cdot&\cdot&\cdot&\cdot&\cdot&\cdot&\cdot&\cdot&\cdot&\cdot&\cdot&\cdot&\cdot&\cdot&\cdot&\cdot&\cdot&\cdot&\cdot&\cdot&\cdot&1&\cdot&-1&\cdot&1\\
&50_x'&\cdot&\cdot&\cdot&\cdot&\cdot&\cdot&\cdot&\cdot&\cdot&\cdot&\cdot&\cdot&\cdot&\cdot&\cdot&\cdot&\cdot&\cdot&\cdot&\cdot&\cdot&\cdot&1&\cdot&\cdot&\cdot\\
&112_z'&\cdot&\cdot&\cdot&\cdot&\cdot&\cdot&\cdot&\cdot&\cdot&\cdot&\cdot&\cdot&\cdot&\cdot&\cdot&\cdot&\cdot&\cdot&\cdot&\cdot&\cdot&\cdot&\cdot&1&-1&-1\\
&35_x'&\cdot&\cdot&\cdot&\cdot&\cdot&\cdot&\cdot&\cdot&\cdot&\cdot&\cdot&\cdot&\cdot&\cdot&\cdot&\cdot&\cdot&\cdot&\cdot&\cdot&\cdot&\cdot&\cdot&\cdot&1&\cdot\\
&1_x'&\cdot&\cdot&\cdot&\cdot&\cdot&\cdot&\cdot&\cdot&\cdot&\cdot&\cdot&\cdot&\cdot&\cdot&\cdot&\cdot&\cdot&\cdot&\cdot&\cdot&\cdot&\cdot&\cdot&\cdot&\cdot&1\\
\end{block}
\end{blockarray}
\]
\end{center}
\end{minipage}
\end{turn}
\restoregeometry
\normalsize
\textbf{\large Finite-dimensional modules, their characters and dimensions:}\\

\begin{align*}
\el(1_x)&=\m(1_x)-\m(112_z)+\m(210_x)-\m(400_z)+\m(525_x)+\m(1400_x)+\m(1050_x)-2\m(2800_z)\\&\qquad-\m(3360_z)+\m(2100_y)+2\m(4536_y)+\m(1400_y)-2\m(2800_z')-\m(3360_z')\\&\qquad+\m(525_x')+\m(1400_x')+\m(1050_x')-\m(400_z')+\m(210_x')-\m(112_z')+\m(1_x')\\
\chi_{\el(1_x)}(t)&=t^{-6}+t^6 + 8(t^{-5}+t^5) + 36(t^{-4}+t^4) + 120(t^{-3}+t^3) + 330(t^{-2}+t^2) + 680(t^{-1}+t) + 1030\\
\dim\el(1_x)&=3380\\
\\
\el(35_x)&=\m(35_x)-\m(112_z)+\m(525_x)+\m(1400_x)-\m(2800_z)-\m(3360_z)+\m(4536_y)+\m(2688_y)\\&\qquad+\m(1400_y)-\m(2800_z')-\m(3360_z')+\m(525_x')+\m(1400_x')-\m(112_z')+\m(35_x')\\
\chi_{\el(35_x)}(t)&=35(t^{-2}+t^2) + 168(t^{-1}+t) + 364\\
\dim\el(35_x)&=770\\
\\
\el(50_x)&=\m(50_x)-\m(400_z)+\m(1400_x)-\m(2800_z)+\m(2100_y)+\m(1400_y)-\m(2800_z')\\&\qquad+\m(1400_x')-\m(400_z')+\m(50_x')\\
\chi_{\el(50_x)}(t)&=50\\
\dim\el(50_x)&=50\\
\end{align*}
\newpage

\subsection{{\Large$\mbf{c=\frac{1}{10}}$}}

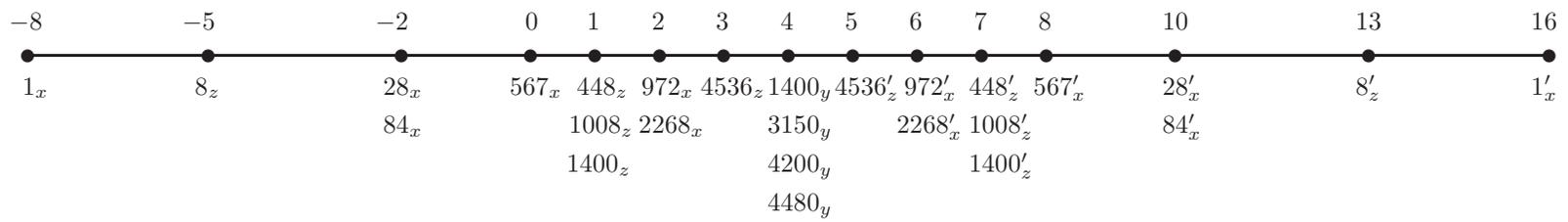
\[
\begin{blockarray}{cccccccccccccccccccccccccccc}
\begin{block}{cc(cccccccccccccccccccccccccc)}
\star&1_x&1&\cdot&1&\cdot&\cdot&\cdot&\cdot&1&\cdot&\cdot&\cdot&\cdot&\cdot&\cdot&\cdot&\cdot&\cdot&\cdot&\cdot&\cdot&\cdot&\cdot&\cdot&\cdot&\cdot&\cdot\\
\star&8_z&\cdot&1&1&\cdot&1&1&\cdot&\cdot&\cdot&\cdot&\cdot&\cdot&\cdot&\cdot&\cdot&\cdot&\cdot&\cdot&\cdot&\cdot&\cdot&\cdot&\cdot&\cdot&\cdot&\cdot\\
\texttt{(2)}&84_x&\cdot&\cdot&1&\cdot&\cdot&1&\cdot&1&\cdot&1&1&\cdot&\cdot&\cdot&\cdot&\cdot&\cdot&\cdot&\cdot&\cdot&\cdot&\cdot&\cdot&\cdot&\cdot&\cdot\\
\star&28_x&\cdot&\cdot&\cdot&1&1&\cdot&1&\cdot&\cdot&\cdot&\cdot&\cdot&\cdot&\cdot&\cdot&\cdot&\cdot&\cdot&\cdot&\cdot&\cdot&\cdot&\cdot&\cdot&\cdot&\cdot\\
\texttt{(1)}&567_x&\cdot&\cdot&\cdot&\cdot&1&1&1&\cdot&1&\cdot&\cdot&\cdot&\cdot&\cdot&\cdot&\cdot&\cdot&\cdot&\cdot&\cdot&\cdot&\cdot&\cdot&\cdot&\cdot&\cdot\\
&1400_z&\cdot&\cdot&\cdot&\cdot&\cdot&1&\cdot&\cdot&1&\cdot&1&1&\cdot&\cdot&\cdot&\cdot&\cdot&\cdot&\cdot&\cdot&\cdot&\cdot&\cdot&\cdot&\cdot&\cdot\\
&1008_z&\cdot&\cdot&\cdot&\cdot&\cdot&\cdot&1&\cdot&1&\cdot&\cdot&\cdot&\cdot&\cdot&1&\cdot&\cdot&\cdot&\cdot&\cdot&\cdot&\cdot&\cdot&\cdot&\cdot&\cdot\\
&448_z&\cdot&\cdot&\cdot&\cdot&\cdot&\cdot&\cdot&1&\cdot&\cdot&1&\cdot&\cdot&1&\cdot&\cdot&\cdot&\cdot&\cdot&\cdot&\cdot&\cdot&\cdot&\cdot&\cdot&\cdot\\
&2268_x&\cdot&\cdot&\cdot&\cdot&\cdot&\cdot&\cdot&\cdot&1&\cdot&\cdot&1&\cdot&\cdot&1&\cdot&1&\cdot&\cdot&\cdot&\cdot&\cdot&\cdot&\cdot&\cdot&\cdot\\
\texttt{(2)}&972_x&\cdot&\cdot&\cdot&\cdot&\cdot&\cdot&\cdot&\cdot&\cdot&1&1&\cdot&1&\cdot&\cdot&\cdot&\cdot&\cdot&\cdot&\cdot&\cdot&\cdot&\cdot&\cdot&\cdot&\cdot\\
&4536_z&\cdot&\cdot&\cdot&\cdot&\cdot&\cdot&\cdot&\cdot&\cdot&\cdot&1&1&1&1&\cdot&1&\cdot&\cdot&\cdot&\cdot&\cdot&\cdot&\cdot&\cdot&\cdot&\cdot\\
&4480_y&\cdot&\cdot&\cdot&\cdot&\cdot&\cdot&\cdot&\cdot&\cdot&\cdot&\cdot&1&\cdot&\cdot&\cdot&1&1&\cdot&1&\cdot&\cdot&\cdot&\cdot&\cdot&\cdot&\cdot\\
&4200_y&\cdot&\cdot&\cdot&\cdot&\cdot&\cdot&\cdot&\cdot&\cdot&\cdot&\cdot&\cdot&1&\cdot&\cdot&1&\cdot&1&\cdot&\cdot&\cdot&\cdot&\cdot&\cdot&\cdot&\cdot\\
&3150_y&\cdot&\cdot&\cdot&\cdot&\cdot&\cdot&\cdot&\cdot&\cdot&\cdot&\cdot&\cdot&\cdot&1&\cdot&1&\cdot&\cdot&\cdot&\cdot&1&\cdot&\cdot&\cdot&\cdot&\cdot\\
&1400_y&\cdot&\cdot&\cdot&\cdot&\cdot&\cdot&\cdot&\cdot&\cdot&\cdot&\cdot&\cdot&\cdot&\cdot&1&\cdot&1&\cdot&\cdot&1&\cdot&\cdot&\cdot&\cdot&\cdot&\cdot\\
&4536_z'&\cdot&\cdot&\cdot&\cdot&\cdot&\cdot&\cdot&\cdot&\cdot&\cdot&\cdot&\cdot&\cdot&\cdot&\cdot&1&\cdot&1&1&\cdot&1&\cdot&1&\cdot&\cdot&\cdot\\
&2268_x'&\cdot&\cdot&\cdot&\cdot&\cdot&\cdot&\cdot&\cdot&\cdot&\cdot&\cdot&\cdot&\cdot&\cdot&\cdot&\cdot&1&\cdot&1&1&\cdot&1&\cdot&\cdot&\cdot&\cdot\\
&972_x'&\cdot&\cdot&\cdot&\cdot&\cdot&\cdot&\cdot&\cdot&\cdot&\cdot&\cdot&\cdot&\cdot&\cdot&\cdot&\cdot&\cdot&1&\cdot&\cdot&\cdot&\cdot&1&\cdot&\cdot&\cdot\\
&1400_z'&\cdot&\cdot&\cdot&\cdot&\cdot&\cdot&\cdot&\cdot&\cdot&\cdot&\cdot&\cdot&\cdot&\cdot&\cdot&\cdot&\cdot&\cdot&1&\cdot&\cdot&1&1&\cdot&1&\cdot\\
&1008_z'&\cdot&\cdot&\cdot&\cdot&\cdot&\cdot&\cdot&\cdot&\cdot&\cdot&\cdot&\cdot&\cdot&\cdot&\cdot&\cdot&\cdot&\cdot&\cdot&1&\cdot&1&\cdot&1&\cdot&\cdot\\
&448_z'&\cdot&\cdot&\cdot&\cdot&\cdot&\cdot&\cdot&\cdot&\cdot&\cdot&\cdot&\cdot&\cdot&\cdot&\cdot&\cdot&\cdot&\cdot&\cdot&\cdot&1&\cdot&1&\cdot&\cdot&1\\
&567_x'&\cdot&\cdot&\cdot&\cdot&\cdot&\cdot&\cdot&\cdot&\cdot&\cdot&\cdot&\cdot&\cdot&\cdot&\cdot&\cdot&\cdot&\cdot&\cdot&\cdot&\cdot&1&\cdot&1&1&\cdot\\
&84_x'&\cdot&\cdot&\cdot&\cdot&\cdot&\cdot&\cdot&\cdot&\cdot&\cdot&\cdot&\cdot&\cdot&\cdot&\cdot&\cdot&\cdot&\cdot&\cdot&\cdot&\cdot&\cdot&1&\cdot&1&1\\
&28_x'&\cdot&\cdot&\cdot&\cdot&\cdot&\cdot&\cdot&\cdot&\cdot&\cdot&\cdot&\cdot&\cdot&\cdot&\cdot&\cdot&\cdot&\cdot&\cdot&\cdot&\cdot&\cdot&\cdot&1&\cdot&\cdot\\
&8_z'&\cdot&\cdot&\cdot&\cdot&\cdot&\cdot&\cdot&\cdot&\cdot&\cdot&\cdot&\cdot&\cdot&\cdot&\cdot&\cdot&\cdot&\cdot&\cdot&\cdot&\cdot&\cdot&\cdot&\cdot&1&\cdot\\
&1_x'&\cdot&\cdot&\cdot&\cdot&\cdot&\cdot&\cdot&\cdot&\cdot&\cdot&\cdot&\cdot&\cdot&\cdot&\cdot&\cdot&\cdot&\cdot&\cdot&\cdot&\cdot&\cdot&\cdot&\cdot&\cdot&1\\
\end{block}
\end{blockarray}
\]

\newgeometry{top=2cm, bottom=2cm}
\begin{turn}{270}
\begin{minipage}{\linewidth}
\small
\begin{center}
\[
\begin{blockarray}{cccccccccccccccccccccccccccc}
\begin{block}{cc(cccccccccccccccccccccccccc)}
\star&1_x&1&\cdot&-1&\cdot&\cdot&1&\cdot&\cdot&-1&1&-1&1&\cdot&1&1&-1&-1&1&1&\cdot&\cdot&\cdot&-1&\cdot&\cdot&1\\
\star&8_z&\cdot&1&-1&\cdot&-1&1&1&1&-1&1&-2&2&1&1&\cdot&-2&-1&1&1&1&1&-1&-1&\cdot&1&\cdot\\
\mathtt{(2)}&84_x&\cdot&\cdot&1&\cdot&\cdot&-1&\cdot&-1&1&-1&2&-2&-1&-1&-1&2&2&-1&-2&-1&-1&1&2&\cdot&-1&-1\\
\star&28_x&\cdot&\cdot&\cdot&1&-1&1&\cdot&\cdot&\cdot&\cdot&-1&\cdot&1&1&\cdot&-1&\cdot&\cdot&1&\cdot&\cdot&-1&\cdot&1&\cdot&\cdot\\
\mathtt{(1)}&567_x&\cdot&\cdot&\cdot&\cdot&1&-1&-1&\cdot&1&\cdot&1&-1&-1&-1&\cdot&2&\cdot&-1&-1&\cdot&-1&1&1&-1&-1&\cdot\\
&1400_z&\cdot&\cdot&\cdot&\cdot&\cdot&1&\cdot&\cdot&-1&\cdot&-1&1&1&1&1&-2&-1&1&2&\cdot&1&-1&-2&1&1&1\\
&1008_z&\cdot&\cdot&\cdot&\cdot&\cdot&\cdot&1&\cdot&-1&\cdot&\cdot&1&\cdot&\cdot&\cdot&-1&\cdot&1&\cdot&\cdot&1&\cdot&-1&\cdot&1&\cdot\\
&448_z&\cdot&\cdot&\cdot&\cdot&\cdot&\cdot&\cdot&1&\cdot&\cdot&-1&1&1&\cdot&\cdot&-1&-1&\cdot&1&1&1&-1&-1&\cdot&1&\cdot\\
&2268_x&\cdot&\cdot&\cdot&\cdot&\cdot&\cdot&\cdot&\cdot&1&\cdot&\cdot&-1&\cdot&\cdot&-1&1&1&-1&-1&\cdot&-1&\cdot&2&\cdot&-1&-1\\
\mathtt{(2)}&972_x&\cdot&\cdot&\cdot&\cdot&\cdot&\cdot&\cdot&\cdot&\cdot&1&-1&1&\cdot&1&\cdot&-1&-1&1&1&1&\cdot&-1&-1&\cdot&1&1\\
&4536_z&\cdot&\cdot&\cdot&\cdot&\cdot&\cdot&\cdot&\cdot&\cdot&\cdot&1&-1&-1&-1&\cdot&2&1&-1&-2&-1&-1&2&2&-1&-2&-1\\
&4480_y&\cdot&\cdot&\cdot&\cdot&\cdot&\cdot&\cdot&\cdot&\cdot&\cdot&\cdot&1&\cdot&\cdot&\cdot&-1&-1&1&1&1&1&-1&-2&\cdot&2&1\\
&4200_y&\cdot&\cdot&\cdot&\cdot&\cdot&\cdot&\cdot&\cdot&\cdot&\cdot&\cdot&\cdot&1&\cdot&\cdot&-1&\cdot&\cdot&1&\cdot&1&-1&-1&1&1&\cdot\\
&3150_y&\cdot&\cdot&\cdot&\cdot&\cdot&\cdot&\cdot&\cdot&\cdot&\cdot&\cdot&\cdot&\cdot&1&\cdot&-1&\cdot&1&1&\cdot&\cdot&-1&-1&1&1&1\\
&1400_y&\cdot&\cdot&\cdot&\cdot&\cdot&\cdot&\cdot&\cdot&\cdot&\cdot&\cdot&\cdot&\cdot&\cdot&1&\cdot&-1&\cdot&1&\cdot&\cdot&\cdot&-1&\cdot&\cdot&1\\
&4536_z'&\cdot&\cdot&\cdot&\cdot&\cdot&\cdot&\cdot&\cdot&\cdot&\cdot&\cdot&\cdot&\cdot&\cdot&\cdot&1&\cdot&-1&-1&\cdot&-1&1&2&-1&-2&-1\\
&2268_x'&\cdot&\cdot&\cdot&\cdot&\cdot&\cdot&\cdot&\cdot&\cdot&\cdot&\cdot&\cdot&\cdot&\cdot&\cdot&\cdot&1&\cdot&-1&-1&\cdot&1&1&\cdot&-1&-1\\
&972_x'&\cdot&\cdot&\cdot&\cdot&\cdot&\cdot&\cdot&\cdot&\cdot&\cdot&\cdot&\cdot&\cdot&\cdot&\cdot&\cdot&\cdot&1&\cdot&\cdot&\cdot&\cdot&-1&\cdot&1&1\\
&1400_z'&\cdot&\cdot&\cdot&\cdot&\cdot&\cdot&\cdot&\cdot&\cdot&\cdot&\cdot&\cdot&\cdot&\cdot&\cdot&\cdot&\cdot&\cdot&1&\cdot&\cdot&-1&-1&1&1&1\\
&1008_z'&\cdot&\cdot&\cdot&\cdot&\cdot&\cdot&\cdot&\cdot&\cdot&\cdot&\cdot&\cdot&\cdot&\cdot&\cdot&\cdot&\cdot&\cdot&\cdot&1&\cdot&-1&\cdot&\cdot&1&\cdot\\
&448_z'&\cdot&\cdot&\cdot&\cdot&\cdot&\cdot&\cdot&\cdot&\cdot&\cdot&\cdot&\cdot&\cdot&\cdot&\cdot&\cdot&\cdot&\cdot&\cdot&\cdot&1&\cdot&-1&\cdot&1&\cdot\\
&567_x'&\cdot&\cdot&\cdot&\cdot&\cdot&\cdot&\cdot&\cdot&\cdot&\cdot&\cdot&\cdot&\cdot&\cdot&\cdot&\cdot&\cdot&\cdot&\cdot&\cdot&\cdot&1&\cdot&-1&-1&\cdot\\
&84_x'&\cdot&\cdot&\cdot&\cdot&\cdot&\cdot&\cdot&\cdot&\cdot&\cdot&\cdot&\cdot&\cdot&\cdot&\cdot&\cdot&\cdot&\cdot&\cdot&\cdot&\cdot&\cdot&1&\cdot&-1&-1\\
&28_x'&\cdot&\cdot&\cdot&\cdot&\cdot&\cdot&\cdot&\cdot&\cdot&\cdot&\cdot&\cdot&\cdot&\cdot&\cdot&\cdot&\cdot&\cdot&\cdot&\cdot&\cdot&\cdot&\cdot&1&\cdot&\cdot\\
&8_z'&\cdot&\cdot&\cdot&\cdot&\cdot&\cdot&\cdot&\cdot&\cdot&\cdot&\cdot&\cdot&\cdot&\cdot&\cdot&\cdot&\cdot&\cdot&\cdot&\cdot&\cdot&\cdot&\cdot&\cdot&1&\cdot\\
&1_x'&\cdot&\cdot&\cdot&\cdot&\cdot&\cdot&\cdot&\cdot&\cdot&\cdot&\cdot&\cdot&\cdot&\cdot&\cdot&\cdot&\cdot&\cdot&\cdot&\cdot&\cdot&\cdot&\cdot&\cdot&\cdot&1\\
\end{block}
\end{blockarray}
\]
\end{center}
\end{minipage}
\end{turn}
\normalsize
\newpage
\restoregeometry

\begin{turn}{270}
\begin{minipage}{\linewidth}
\small
\begin{center}
\begin{picture}(600,100)
\put(5,80){\line(1,0){590}}
\put(5,80){\circle*{5}}
\put(75,80){\circle*{5}}
\put(150,80){\circle*{5}}
\put(200,80){\circle*{5}}
\put(225,80){\circle*{5}}
\put(250,80){\circle*{5}}
\put(275,80){\circle*{5}}
\put(300,80){\circle*{5}}
\put(325,80){\circle*{5}}
\put(350,80){\circle*{5}}
\put(375,80){\circle*{5}}
\put(400,80){\circle*{5}}
\put(450,80){\circle*{5}}
\put(525,80){\circle*{5}}
\put(595,80){\circle*{5}}
\put(-2,90){$-8$}
\put(65,90){$-5$}
\put(140,90){$-2$}
\put(198,90){$0$}
\put(222,90){$1$}
\put(247,90){$2$}
\put(272,90){$3$}
\put(297,90){$4$}
\put(322,90){$5$}
\put(347,90){$6$}
\put(372,90){$7$}
\put(397,90){$8$}
\put(445,90){$10$}
\put(520,90){$13$}
\put(588,90){$16$}
\put(3,65){$1_x$}
\put(70,65){$8_z$}
\put(143,65){$28_x$}
\put(143,50){$84_x$}
\put(192,65){$567_x$}
\put(218,65){$448_z$}
\put(215,50){$1008_z$}
\put(214,35){$1400_z$}
\put(243,65){$972_x$}
\put(242,50){$2268_x$}
\put(266,65){$4536_z$}
\put(292,65){$1400_y$}
\put(292,50){$3150_y$}
\put(292,35){$4200_y$}
\put(292,20){$4480_y$}
\put(318,65){$4536_z'$}
\put(345,65){$972_x'$}
\put(342,50){$2268_x'$}
\put(370,65){$448_z'$}
\put(370,50){$1008_z'$}
\put(370,35){$1400_z'$}
\put(395,65){$567_x'$}
\put(445,65){$28_x'$}
\put(445,50){$84_x'$}
\put(520,65){$8_z'$}
\put(588,65){$1_x'$}
\end{picture}
\end{center}
\end{minipage}
\end{turn}

\restoregeometry

\textbf{\large Finite-dimensional modules, their characters and dimensions:}\\
\begin{align*}
&\el(1_x)=\m(1_x)-\m(84_x)+\m(1400_z)-\m(2268_x)+\m(972_x)-\m(4536_z)+\m(4480_y)+\m(3150_y)\\&\qquad+\m(1400_y)-\m(4536_z')-\m(2268_x')+\m(972_x')+\m(1400_z')-\m(84_x')+\m(1_x')\\
&\chi_{\el(1_x)}(t)=t^{-8}+t^8 + 8(t^{-7}+t^7) + 36(t^{-6}+t^6) + 120(t^{-5}+t^5) + 330(t^{-4}+t^4) + 792(t^{-3}+t^3)\\&\qquad + 1632(t^{-2}+t^2) + 2760(t^{-1}+t) + 3411\\
&\dim\el(1_x)=14,769\\
\\
&\el(8_z)=\m(8_z)-\m(84_x)-\m(567_x)+\m(1400_z)+\m(1008_z)+\m(448_z)-\m(2268_x)+\m(972_x)\\&\qquad-2\m(4536_z)+2\m(4480_y)+\m(4200_y)+\m(3150_y)-2\m(4536_z')-\m(2268_x')+\m(972_x')\\
&\qquad+\m(1400_z')+\m(1008_z')+\m(448_z')-\m(567_x')-\m(84_x')+\m(8_z')\\
&\chi_{\el(8_z)}(t)=8(t^{-5}+t^5) + 64(t^{-4}+t^4) + 288(t^{-3}+t^3) + 876(t^{-2}+t^2) + 1968(t^{-1}+t) + 2745\\
&\dim\el(8_z)=9153\\
\\
&\el(28_x)=\m(28_x)-\m(567_x)+\m(1400_z)-\m(4536_z)+\m(4200_y)+\m(3150_y)-\m(4536_z')\\&\qquad+\m(1400_z')-\m(567_x')+\m(28_x')\\
&\chi_{\el(28_x)}(t)=28(t^{-2}+t^2) + 224(t^{-1}+t) + 441\\
&\dim\el(28_x)=945\\
\end{align*}

\newgeometry{left=0cm,right=0cm,top=4cm}
\subsection{{\Large$\mbf{c=\frac{1}{8}}$}}
\small
\begin{center}
\vspace*{1cm}
\begin{picture}(600,100)
\put(5,80){\line(1,0){590}}
\put(5,80){\circle*{5}}
\put(75,80){\circle*{5}}
\put(150,80){\circle*{5}}
\put(175,80){\circle*{5}}
\put(225,80){\circle*{5}}
\put(275,80){\circle*{5}}
\put(300,80){\circle*{5}}
\put(325,80){\circle*{5}}
\put(375,80){\circle*{5}}
\put(425,80){\circle*{5}}
\put(450,80){\circle*{5}}
\put(525,80){\circle*{5}}
\put(595,80){\circle*{5}}
\put(-2,90){$-11$}
\put(65,90){$-5$}
\put(140,90){$-2$}
\put(172,90){$-1$}
\put(222,90){$1$}
\put(272,90){$3$}
\put(297,90){$4$}
\put(322,90){$5$}
\put(372,90){$7$}
\put(422,90){$9$}
\put(445,90){$10$}
\put(520,90){$13$}
\put(588,90){$19$}
\put(3,65){$1_x$}
\put(70,65){$35_x$}
\put(140,65){$160_z$}
\put(172,65){$567_x$}
\put(216,65){$525_x$}
\put(215,50){$1575_x$}
\put(214,35){$1400_x$}
\put(216,20){$175_x$}
\put(266,65){$6075_x$}
\put(266,50){$2835_x$}
\put(292,65){$7168_w$}
\put(292,50){$5600_w$}
\put(292,35){$2016_w$}
\put(318,65){$6075_x'$}
\put(318,50){$2835_x'$}
\put(370,65){$525_x'$}
\put(370,50){$1575_x'$}
\put(370,35){$1400_x'$}
\put(370,20){$175_x'$}
\put(420,65){$567_x'$}
\put(445,50){$160_z'$}
\put(520,65){$35_x'$}
\put(588,65){$1_x'$}
\end{picture}
\end{center}

\small
\[
\begin{blockarray}{ccccccccccccccccccccccccc}
\begin{block}{cc(ccccccccccccccccccccccc)}
\star&1_x&1&1&\cdot&\cdot&\cdot&\cdot&\cdot&1&\cdot&1&\cdot&\cdot&\cdot&\cdot&\cdot&\cdot&\cdot&\cdot&\cdot&\cdot&\cdot&\cdot&\cdot\\
\texttt{(3)}&35_x&\cdot&1&\cdot&1&\cdot&\cdot&1&\cdot&\cdot&1&\cdot&\cdot&\cdot&\cdot&\cdot&\cdot&\cdot&\cdot&\cdot&\cdot&\cdot&\cdot&\cdot\\
\star&160_z&\cdot&\cdot&1&1&\cdot&1&\cdot&\cdot&\cdot&\cdot&\cdot&\cdot&\cdot&\cdot&\cdot&\cdot&\cdot&\cdot&\cdot&\cdot&\cdot&\cdot&\cdot\\
\texttt{(3)}&567_x&\cdot&\cdot&\cdot&1&1&1&1&\cdot&1&\cdot&\cdot&\cdot&\cdot&\cdot&\cdot&\cdot&\cdot&\cdot&\cdot&\cdot&\cdot&\cdot&\cdot\\
\texttt{(3)}&525_x&\cdot&\cdot&\cdot&\cdot&1&\cdot&\cdot&\cdot&1&\cdot&\cdot&\cdot&\cdot&\cdot&\cdot&1&\cdot&\cdot&\cdot&\cdot&\cdot&\cdot&\cdot\\
&1575_x&\cdot&\cdot&\cdot&\cdot&\cdot&1&\cdot&\cdot&1&\cdot&\cdot&1&\cdot&\cdot&\cdot&\cdot&\cdot&\cdot&\cdot&\cdot&\cdot&\cdot&\cdot\\
&1400_x&\cdot&\cdot&\cdot&\cdot&\cdot&\cdot&1&\cdot&1&1&1&\cdot&\cdot&\cdot&\cdot&\cdot&\cdot&\cdot&\cdot&\cdot&\cdot&\cdot&\cdot\\
\texttt{(1)}&175_x&\cdot&\cdot&\cdot&\cdot&\cdot&\cdot&\cdot&1&\cdot&1&\cdot&\cdot&1&\cdot&\cdot&\cdot&\cdot&\cdot&\cdot&\cdot&\cdot&\cdot&\cdot\\
&6075_x&\cdot&\cdot&\cdot&\cdot&\cdot&\cdot&\cdot&\cdot&1&\cdot&1&1&\cdot&1&\cdot&1&\cdot&\cdot&\cdot&\cdot&\cdot&\cdot&\cdot\\
&2835_x&\cdot&\cdot&\cdot&\cdot&\cdot&\cdot&\cdot&\cdot&\cdot&1&1&\cdot&1&\cdot&1&\cdot&\cdot&\cdot&\cdot&\cdot&\cdot&\cdot&\cdot\\
&7168_w&\cdot&\cdot&\cdot&\cdot&\cdot&\cdot&\cdot&\cdot&\cdot&\cdot&1&\cdot&\cdot&1&1&\cdot&\cdot&1&\cdot&\cdot&\cdot&\cdot&\cdot\\
&5600_w&\cdot&\cdot&\cdot&\cdot&\cdot&\cdot&\cdot&\cdot&\cdot&\cdot&\cdot&1&\cdot&1&\cdot&\cdot&1&\cdot&\cdot&\cdot&\cdot&\cdot&\cdot\\
&2016_w&\cdot&\cdot&\cdot&\cdot&\cdot&\cdot&\cdot&\cdot&\cdot&\cdot&\cdot&\cdot&1&\cdot&1&\cdot&\cdot&\cdot&1&\cdot&\cdot&\cdot&\cdot\\
&6075_x'&\cdot&\cdot&\cdot&\cdot&\cdot&\cdot&\cdot&\cdot&\cdot&\cdot&\cdot&\cdot&\cdot&1&\cdot&1&1&1&\cdot&1&\cdot&\cdot&\cdot\\
&2835_x'&\cdot&\cdot&\cdot&\cdot&\cdot&\cdot&\cdot&\cdot&\cdot&\cdot&\cdot&\cdot&\cdot&\cdot&1&\cdot&\cdot&1&1&\cdot&\cdot&\cdot&1\\
&525_x'&\cdot&\cdot&\cdot&\cdot&\cdot&\cdot&\cdot&\cdot&\cdot&\cdot&\cdot&\cdot&\cdot&\cdot&\cdot&1&\cdot&\cdot&\cdot&1&\cdot&\cdot&\cdot\\
&1575_x'&\cdot&\cdot&\cdot&\cdot&\cdot&\cdot&\cdot&\cdot&\cdot&\cdot&\cdot&\cdot&\cdot&\cdot&\cdot&\cdot&1&\cdot&\cdot&1&1&\cdot&\cdot\\
&1400_x'&\cdot&\cdot&\cdot&\cdot&\cdot&\cdot&\cdot&\cdot&\cdot&\cdot&\cdot&\cdot&\cdot&\cdot&\cdot&\cdot&\cdot&1&\cdot&1&\cdot&1&1\\
&175_x'&\cdot&\cdot&\cdot&\cdot&\cdot&\cdot&\cdot&\cdot&\cdot&\cdot&\cdot&\cdot&\cdot&\cdot&\cdot&\cdot&\cdot&\cdot&1&\cdot&\cdot&\cdot&1\\
&567_x'&\cdot&\cdot&\cdot&\cdot&\cdot&\cdot&\cdot&\cdot&\cdot&\cdot&\cdot&\cdot&\cdot&\cdot&\cdot&\cdot&\cdot&\cdot&\cdot&1&1&1&\cdot\\
&160_z'&\cdot&\cdot&\cdot&\cdot&\cdot&\cdot&\cdot&\cdot&\cdot&\cdot&\cdot&\cdot&\cdot&\cdot&\cdot&\cdot&\cdot&\cdot&\cdot&\cdot&1&\cdot&\cdot\\
&35_x'&\cdot&\cdot&\cdot&\cdot&\cdot&\cdot&\cdot&\cdot&\cdot&\cdot&\cdot&\cdot&\cdot&\cdot&\cdot&\cdot&\cdot&\cdot&\cdot&\cdot&\cdot&1&1\\
&1_x'&\cdot&\cdot&\cdot&\cdot&\cdot&\cdot&\cdot&\cdot&\cdot&\cdot&\cdot&\cdot&\cdot&\cdot&\cdot&\cdot&\cdot&\cdot&\cdot&\cdot&\cdot&\cdot&1\\
\end{block}
\end{blockarray}
\]

\newgeometry{top=3cm}
\begin{turn}{270}
\begin{minipage}{\linewidth}
\small
\[
\begin{blockarray}{ccccccccccccccccccccccccc}
\begin{block}{cc(ccccccccccccccccccccccc)}
\star&1_x&1&-1&\cdot&1&-1&-1&\cdot&-1&1&1&-2&\cdot&\cdot&1&1&-1&-1&\cdot&-1&1&\cdot&-1&1\\
\texttt{(3)}&35_x&\cdot&1&\cdot&-1&1&1&\cdot&\cdot&-1&-1&2&\cdot&1&-1&-2&1&1&1&1&-2&1&1&-1\\
\star&160_z&\cdot&\cdot&1&-1&1&\cdot&1&\cdot&-1&-1&1&1&1&-1&-1&1&\cdot&1&\cdot&-1&1&\cdot&\cdot\\
\texttt{(3)}&567_x&\cdot&\cdot&\cdot&1&-1&-1&-1&\cdot&2&1&-2&-1&-1&1&2&-2&\cdot&-1&-1&2&-2&-1&1\\
\texttt{(3)}&525_x&\cdot&\cdot&\cdot&\cdot&1&\cdot&\cdot&\cdot&-1&\cdot&1&1&\cdot&-1&-1&1&\cdot&1&1&-1&1&\cdot&-1\\
&1575_x&\cdot&\cdot&\cdot&\cdot&\cdot&1&\cdot&\cdot&-1&\cdot&1&\cdot&\cdot&\cdot&-1&1&\cdot&\cdot&1&-1&1&1&-1\\
&1400_x&\cdot&\cdot&\cdot&\cdot&\cdot&\cdot&1&\cdot&-1&-1&1&1&1&-1&-1&2&\cdot&1&\cdot&-2&2&1&-1\\
\texttt{(1)}&175_x&\cdot&\cdot&\cdot&\cdot&\cdot&\cdot&\cdot&1&\cdot&-1&1&\cdot&\cdot&-1&\cdot&1&1&\cdot&\cdot&-1&\cdot&1&-1\\
&6075_x&\cdot&\cdot&\cdot&\cdot&\cdot&\cdot&\cdot&\cdot&1&\cdot&-1&-1&\cdot&1&1&-2&\cdot&-1&-1&2&-2&-1&2\\
&2835_x&\cdot&\cdot&\cdot&\cdot&\cdot&\cdot&\cdot&\cdot&\cdot&1&-1&\cdot&-1&1&1&-1&-1&-1&\cdot&2&-1&-1&1\\
&7168_w&\cdot&\cdot&\cdot&\cdot&\cdot&\cdot&\cdot&\cdot&\cdot&\cdot&1&\cdot&\cdot&-1&-1&1&1&1&1&-2&1&1&-2\\
&5600_w&\cdot&\cdot&\cdot&\cdot&\cdot&\cdot&\cdot&\cdot&\cdot&\cdot&\cdot&1&\cdot&-1&\cdot&1&\cdot&1&\cdot&-1&1&\cdot&-1\\
&2016_w&\cdot&\cdot&\cdot&\cdot&\cdot&\cdot&\cdot&\cdot&\cdot&\cdot&\cdot&\cdot&1&\cdot&-1&\cdot&\cdot&1&\cdot&-1&1&\cdot&\cdot\\
&6075_x'&\cdot&\cdot&\cdot&\cdot&\cdot&\cdot&\cdot&\cdot&\cdot&\cdot&\cdot&\cdot&\cdot&1&\cdot&-1&-1&-1&\cdot&2&-1&-1&2\\
&2835_x'&\cdot&\cdot&\cdot&\cdot&\cdot&\cdot&\cdot&\cdot&\cdot&\cdot&\cdot&\cdot&\cdot&\cdot&1&\cdot&\cdot&-1&-1&1&-1&\cdot&1\\
&525_x'&\cdot&\cdot&\cdot&\cdot&\cdot&\cdot&\cdot&\cdot&\cdot&\cdot&\cdot&\cdot&\cdot&\cdot&\cdot&1&\cdot&\cdot&\cdot&-1&1&1&-1\\
&1575_x'&\cdot&\cdot&\cdot&\cdot&\cdot&\cdot&\cdot&\cdot&\cdot&\cdot&\cdot&\cdot&\cdot&\cdot&\cdot&\cdot&1&\cdot&\cdot&-1&\cdot&1&-1\\
&1400_x'&\cdot&\cdot&\cdot&\cdot&\cdot&\cdot&\cdot&\cdot&\cdot&\cdot&\cdot&\cdot&\cdot&\cdot&\cdot&\cdot&\cdot&1&\cdot&-1&1&\cdot&-1\\
&175_x'&\cdot&\cdot&\cdot&\cdot&\cdot&\cdot&\cdot&\cdot&\cdot&\cdot&\cdot&\cdot&\cdot&\cdot&\cdot&\cdot&\cdot&\cdot&1&\cdot&\cdot&\cdot&-1\\
&567_x'&\cdot&\cdot&\cdot&\cdot&\cdot&\cdot&\cdot&\cdot&\cdot&\cdot&\cdot&\cdot&\cdot&\cdot&\cdot&\cdot&\cdot&\cdot&\cdot&1&-1&-1&1\\
&160_z'&\cdot&\cdot&\cdot&\cdot&\cdot&\cdot&\cdot&\cdot&\cdot&\cdot&\cdot&\cdot&\cdot&\cdot&\cdot&\cdot&\cdot&\cdot&\cdot&\cdot&1&\cdot&\cdot\\
&35_x'&\cdot&\cdot&\cdot&\cdot&\cdot&\cdot&\cdot&\cdot&\cdot&\cdot&\cdot&\cdot&\cdot&\cdot&\cdot&\cdot&\cdot&\cdot&\cdot&\cdot&\cdot&1&-1\\
&1_x'&\cdot&\cdot&\cdot&\cdot&\cdot&\cdot&\cdot&\cdot&\cdot&\cdot&\cdot&\cdot&\cdot&\cdot&\cdot&\cdot&\cdot&\cdot&\cdot&\cdot&\cdot&\cdot&1\\
\end{block}
\end{blockarray}
\]
\end{minipage}
\end{turn}
\restoregeometry
\normalsize

\textbf{\large Finite-dimensional modules, their characters and dimensions:}\\
\begin{align*}
\el(1_x)&=\m(1_x)-\m(35_x)+\m(567_x)-\m(525_x)-\m(1575_x)-\m(175_x)+\m(6075_x)+\m(2835_x)\\&\qquad-2\m(7168_w)+\m(6075_x')+\m(2835_x')-\m(525_x')-\m(1575_x')-\m(175_x')\\&\qquad+\m(567_x')-\m(35_x')+\m(1_x')\\
\chi_{\el(1_x)}(t)&=t^{-11}+t^{11} + 8(t^{-10}+t^{10}) + 36(t^{-9}+t^9) + 120(t^{-8}+t^8) + 330(t^{-7}+t^7) + 792(t^{-6}+t^6) \\&\qquad+ 1681(t^{-5}+t^5) + 3152(t^{-4}+t^4) + 5175(t^{-3}+t^3) + 7240(t^{-2}+t^2) + 8465(t^{-1}+t) + 8640\\
\dim\el(1_x)&=62,640\\
\\
\el(160_z)&=\m(160_z)-\m(567_x)+\m(525_x)+\m(1400_x)-\m(6075_x)-\m(2835_x)+\m(7168_w)+\m(5600_w)\\&\qquad+\m(2016_w)-\m(6075_x')-\m(2835_x')+\m(525_x')+\m(1400_x')-\m(567_x')+\m(160_z')\\
\chi_{\el(160_z)}(t)&=160(t^{-2}+t^2) + 713(t^{-1}+t) + 1224\\
\dim\el(160_z)&=2970\\
\\
\end{align*}

\subsection{{\Large$\mbf{c=\frac{1}{6}}$}}

There is a large principal block containing the finite-dimensional $\el(\mathrm{Triv})$. The principal block is not studied in this paper. There is also a smaller block whose minimal support irreducibles have support of dimension $2$, and its decomposition matrix is given below.

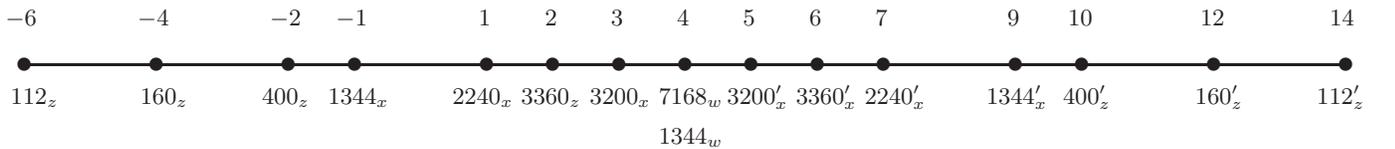
\footnotesize
\begin{center}
\begin{picture}(500,65)
\put(0,45){\line(1,0){500}}
\put(0,45){\circle*{5}}
\put(50,45){\circle*{5}}
\put(100,45){\circle*{5}}
\put(125,45){\circle*{5}}
\put(175,45){\circle*{5}}
\put(200,45){\circle*{5}}
\put(225,45){\circle*{5}}
\put(250,45){\circle*{5}}
\put(275,45){\circle*{5}}
\put(300,45){\circle*{5}}
\put(325,45){\circle*{5}}
\put(375,45){\circle*{5}}
\put(400,45){\circle*{5}}
\put(450,45){\circle*{5}}
\put(500,45){\circle*{5}}
\put(-7,60){$-6$}
\put(43,60){$-4$}
\put(93,60){$-2$}
\put(118,60){$-1$}
\put(172,60){$1$}
\put(197,60){$2$}
\put(222,60){$3$}
\put(247,60){$4$}
\put(272,60){$5$}
\put(297,60){$6$}
\put(322,60){$7$}
\put(372,60){$9$}
\put(395,60){$10$}
\put(445,60){$12$}
\put(494,60){$14$}
\put(-5,30){$112_z$}
\put(44,30){$160_z$}
\put(90,30){$400_z$}
\put(115,30){$1344_x$}
\put(162,30){$2240_x$}
\put(188,30){$3360_z$}
\put(214,30){$3200_x$}
\put(240,30){$7168_w$}
\put(240,15){$1344_w$}
\put(266,30){$3200_x'$}
\put(292,30){$3360_x'$}
\put(318,30){$2240_x'$}
\put(364,30){$1344_x'$}
\put(393,30){$400_z'$}
\put(443,30){$160_z'$}
\put(489,30){$112_z'$}
\end{picture}
\end{center}
\normalsize

\[
\begin{blockarray}{cccccccccccccccccc}
\begin{block}{cc(cccccccccccccccc)}
\mathtt{(2)}&112_z&1&\cdot&1&1&1&\cdot&\cdot&\cdot&\cdot&\cdot&\cdot&\cdot&\cdot&\cdot&\cdot&\cdot\\
\mathtt{(2)}&160_z&\cdot&1&\cdot&1&\cdot&1&\cdot&\cdot&\cdot&\cdot&\cdot&\cdot&\cdot&\cdot&\cdot&\cdot\\
\mathtt{(3)}&400_z&\cdot&\cdot&1&\cdot&1&\cdot&\cdot&\cdot&1&\cdot&\cdot&\cdot&\cdot&\cdot&\cdot&\cdot\\
\mathtt{(4)}&1344_x&\cdot&\cdot&\cdot&1&1&1&1&1&\cdot&\cdot&\cdot&\cdot&\cdot&\cdot&\cdot&\cdot\\
&2240_x&\cdot&\cdot&\cdot&\cdot&1&\cdot&\cdot&1&1&\cdot&\cdot&1&\cdot&\cdot&\cdot&\cdot\\
&3360_z&\cdot&\cdot&\cdot&\cdot&\cdot&1&\cdot&1&\cdot&1&\cdot&\cdot&\cdot&\cdot&\cdot&\cdot\\
\mathtt{(4)}&3200_x&\cdot&\cdot&\cdot&\cdot&\cdot&\cdot&1&1&\cdot&\cdot&1&\cdot&\cdot&\cdot&\cdot&\cdot\\
&7168_w&\cdot&\cdot&\cdot&\cdot&\cdot&\cdot&\cdot&1&\cdot&1&1&1&1&\cdot&\cdot&\cdot\\
&1344_w&\cdot&\cdot&\cdot&\cdot&\cdot&\cdot&\cdot&\cdot&1&\cdot&\cdot&1&\cdot&1&\cdot&\cdot\\
&3200_x'&\cdot&\cdot&\cdot&\cdot&\cdot&\cdot&\cdot&\cdot&\cdot&1&\cdot&\cdot&1&\cdot&\cdot&\cdot\\
&3360_z'&\cdot&\cdot&\cdot&\cdot&\cdot&\cdot&\cdot&\cdot&\cdot&\cdot&1&\cdot&1&\cdot&1&\cdot\\
&2240_x'&\cdot&\cdot&\cdot&\cdot&\cdot&\cdot&\cdot&\cdot&\cdot&\cdot&\cdot&1&1&1&\cdot&1\\
&1344_x'&\cdot&\cdot&\cdot&\cdot&\cdot&\cdot&\cdot&\cdot&\cdot&\cdot&\cdot&\cdot&1&\cdot&1&1\\
&400_z'&\cdot&\cdot&\cdot&\cdot&\cdot&\cdot&\cdot&\cdot&\cdot&\cdot&\cdot&\cdot&\cdot&1&\cdot&1\\
&160_z'&\cdot&\cdot&\cdot&\cdot&\cdot&\cdot&\cdot&\cdot&\cdot&\cdot&\cdot&\cdot&\cdot&\cdot&1&\cdot\\
&112_z'&\cdot&\cdot&\cdot&\cdot&\cdot&\cdot&\cdot&\cdot&\cdot&\cdot&\cdot&\cdot&\cdot&\cdot&\cdot&1\\
\end{block}\\
\begin{block}{cc(cccccccccccccccc)}
\mathtt{(2)}&112_z&1&\cdot&-1&-1&1&1&1&-2&\cdot&1&1&1&-1&-1&\cdot&1\\
\mathtt{(2)}&160_z&\cdot&1&\cdot&-1&1&\cdot&1&-1&-1&1&\cdot&1&-1&\cdot&1&\cdot\\
\mathtt{(3)}&400_z&\cdot&\cdot&1&\cdot&-1&\cdot&\cdot&1&\cdot&-1&-1&\cdot&1&\cdot&\cdot&-1\\
\mathtt{(4)}&1344_x&\cdot&\cdot&\cdot&1&-1&-1&-1&2&1&-1&-1&-2&2&1&-1&-1\\
&2240_x&\cdot&\cdot&\cdot&\cdot&1&\cdot&\cdot&-1&-1&1&1&1&-2&\cdot&1&1\\
&3360_z&\cdot&\cdot&\cdot&\cdot&\cdot&1&\cdot&-1&\cdot&\cdot&1&1&-1&-1&\cdot&1\\
\mathtt{(4)}&3200_x&\cdot&\cdot&\cdot&\cdot&\cdot&\cdot&1&-1&\cdot&1&\cdot&1&-1&-1&1&1\\
&7168_w&\cdot&\cdot&\cdot&\cdot&\cdot&\cdot&\cdot&1&\cdot&-1&-1&-1&2&1&-1&-2\\
&1344_w&\cdot&\cdot&\cdot&\cdot&\cdot&\cdot&\cdot&\cdot&1&\cdot&\cdot&-1&1&\cdot&-1&\cdot\\
&3200_x'&\cdot&\cdot&\cdot&\cdot&\cdot&\cdot&\cdot&\cdot&\cdot&1&\cdot&\cdot&-1&\cdot&1&1\\
&3360_z'&\cdot&\cdot&\cdot&\cdot&\cdot&\cdot&\cdot&\cdot&\cdot&\cdot&1&\cdot&-1&\cdot&\cdot&1\\
&2240_x'&\cdot&\cdot&\cdot&\cdot&\cdot&\cdot&\cdot&\cdot&\cdot&\cdot&\cdot&1&-1&-1&1&1\\
&1344_x'&\cdot&\cdot&\cdot&\cdot&\cdot&\cdot&\cdot&\cdot&\cdot&\cdot&\cdot&\cdot&1&\cdot&-1&-1\\
&400_z'&\cdot&\cdot&\cdot&\cdot&\cdot&\cdot&\cdot&\cdot&\cdot&\cdot&\cdot&\cdot&\cdot&1&\cdot&-1\\
&160_z'&\cdot&\cdot&\cdot&\cdot&\cdot&\cdot&\cdot&\cdot&\cdot&\cdot&\cdot&\cdot&\cdot&\cdot&1&\cdot\\
&112_z'&\cdot&\cdot&\cdot&\cdot&\cdot&\cdot&\cdot&\cdot&\cdot&\cdot&\cdot&\cdot&\cdot&\cdot&\cdot&1\\
\end{block}\\
\end{blockarray}
\]

All missing entries can be recovered from the Hecke algebra decomposition matrix using Lemmas \ref{RR} and \ref{dim Hom} except for $[\m(112_z):\el(3200_x)]$ and $[\m(160_z):\el(3200_x)]$. Let 
$$\alpha:=[\m(160_z):\el(3200_x)],\qquad\beta:=[\m(112_z):\el(3200_x)]$$
and write the decompositions of $\el(160_z)$ and $\el(112_z)$ into Vermas in terms of $\alpha$ and $\beta$:
\begin{align*}
\el(160_z)&=\m(160_z)-\m(1344_x)+\m(2240_x)+(1-\alpha)\m(3200_x)+(\alpha-1)\m(7168_w)-\m(1344_w)\\
&\qquad+\m(3200_x')+(1-\alpha)\m(2240_x')-\m(1344_x')+\alpha\m(400_z')+\m(160_z')\\
\\
\el(112_z)&=\m(112_z)-\m(400_z)-\m(1344_x)+\m(2240_x)+\m(3360_z)+(1-\beta)\m(3200_x)\\
&\qquad+(\beta-2)\m(7168_w)+(1-\beta)\m(3200_x')+\m(3360_z')+(1-\beta)\m(2240_x')\\
&\qquad+(\beta-1)\m(1344_x')+(\beta-1)\m(400_z')-\beta\m(160_z')+(1-\beta)\m(112_z')\\
\end{align*}
In the first matrix we see that $[\m(3200_x'):\el(\sigma)]=0$ unless $\sigma=1344_x'$, and so\\ $\dim\Hom(\m(\sigma),\m(3200_x'))=0$ unless $\sigma=1344_x'$ or possibly $\sigma$ is the lowest weight of a simple appearing in the decomposition of $\m(1344_x')$. However, from the decompositions of $\el(3200_x')=\m(3200_x')-\m(1344_x')+\m(160_z')+\m(112_z')$ and $\el(1344_x')=\m(1344_x')-\m(160_z')-\m(112_z')$, it must be that $$\dim\Hom(\m(160_z'),\m(1344_x')=\dim\Hom(\m(112_z'),\m(1344_x'))=0$$ as well. So $\dim\Hom(\m(\sigma),\m(3200_x'))=0$ unless $\sigma=1344_x'$. By Lemma \ref{dim Hom}, $\alpha$ and $\beta$ must be either $0$ or $1$.

By considering dimensions of supports we can show that $\alpha$ and $\beta$ must both be $0$. If $\alpha=1$, the character of $\el(160_z)$ would have a pole of order $8$, i.e. $\el(160_z)$ would have full support. But this is impossible since $\el(160_z)$ is killed by KZ functor. Therefore $\alpha=0$. If $\beta=1$ then the character of $\el(112_z)$ would have a pole of order $4$, i.e. the dimension of support of $\el(112_z)$ would be $4$-dimensional. But consider $$\Res^{E_8}_{E_7} 112_z=1_a\oplus27_a\oplus7_a'\oplus21_b'\oplus56a'$$
The dimension of support of the simple $H_{\frac{1}{6}}(E_7)$-module $\el(56_a')$ is $2$, and therefore $\Ind_{E_7}^{E_8}\el(56_a')$ has $3$-dimensional support. But 
$$\Ind\el(56_a')=\Ind\m(56_a')-\Ind\m(120_a)+\Ind\m(336_a')-\Ind\m(336_a)+\Ind\m(120_a')-\Ind\m(56_a)$$ and from the formula for $\Res 112_z$ above, $\m(112_z)$ appears exactly once in the formula for $\Ind\el(56_a')$, in $\Ind\m(56_a')$. Since $112_z$ is leftmost in its block with respect to the $h_c$-ordering, it follows that $\el(112_z)$ appears as a summand of $\Ind\el(56_a')$ in the Grothendieck group. Since $\dim\Supp\Ind\el(56_a')=3$, $\dim\Supp\el(112_z)\leq3$, and so $\beta=0$.

\subsection{{\Large$\mbf{c=\frac{1}{5}}$}} There are two nontrivial blocks, one containing $\el(\mathrm{Triv})=\el(1_x)$ and the other containing $\el(8_z)$, where $8_z$ denotes the standard representation of the Weyl group $E_8$. The blocks appear identical after relabeling: each block contains twenty irreducible representations, of which two are finite-dimensional and four have $4$-dimensional support. 

\subsubsection{The principal block}

\small
\[
\begin{blockarray}{cccccccccccccccccccccc}
\begin{block}{cc(cccccccccccccccccccc)}
\star&1_x&1&1&\cdot&\cdot&\cdot&\cdot&\cdot&\cdot&\cdot&\cdot&\cdot&1&\cdot&\cdot&\cdot&\cdot&\cdot&\cdot&\cdot&\cdot\\
\mathtt{(4)}&84_x&\cdot&1&\cdot&\cdot&1&\cdot&1&\cdot&\cdot&\cdot&\cdot&1&\cdot&\cdot&\cdot&\cdot&\cdot&\cdot&\cdot&\cdot\\
\star&28_x&\cdot&\cdot&1&1&\cdot&\cdot&\cdot&\cdot&\cdot&\cdot&1&\cdot&\cdot&\cdot&\cdot&\cdot&\cdot&\cdot&\cdot&\cdot\\
\mathtt{(4)}&567_x&\cdot&\cdot&\cdot&1&1&1&\cdot&1&\cdot&\cdot&1&\cdot&\cdot&\cdot&\cdot&\cdot&\cdot&\cdot&\cdot&\cdot\\
\mathtt{(4)}&1344_x&\cdot&\cdot&\cdot&\cdot&1&\cdot&1&1&\cdot&1&\cdot&\cdot&\cdot&\cdot&\cdot&\cdot&\cdot&\cdot&\cdot&\cdot\\
\mathtt{(4)}&2268_x&\cdot&\cdot&\cdot&\cdot&\cdot&1&\cdot&1&1&\cdot&\cdot&\cdot&\cdot&\cdot&\cdot&\cdot&\cdot&\cdot&\cdot&\cdot\\
&972_x&\cdot&\cdot&\cdot&\cdot&\cdot&\cdot&1&\cdot&\cdot&1&\cdot&1&\cdot&\cdot&1&\cdot&\cdot&\cdot&\cdot&\cdot\\
&4096_z&\cdot&\cdot&\cdot&\cdot&\cdot&\cdot&\cdot&1&1&1&1&\cdot&1&\cdot&\cdot&\cdot&\cdot&\cdot&\cdot&\cdot\\
&4536_y&\cdot&\cdot&\cdot&\cdot&\cdot&\cdot&\cdot&\cdot&1&\cdot&\cdot&\cdot&1&1&\cdot&\cdot&\cdot&\cdot&\cdot&\cdot\\
&2688_y&\cdot&\cdot&\cdot&\cdot&\cdot&\cdot&\cdot&\cdot&\cdot&1&\cdot&\cdot&1&\cdot&1&1&\cdot&\cdot&\cdot&\cdot\\
&1134_y&\cdot&\cdot&\cdot&\cdot&\cdot&\cdot&\cdot&\cdot&\cdot&\cdot&1&\cdot&1&\cdot&\cdot&\cdot&\cdot&\cdot&1&\cdot\\
&168_y&\cdot&\cdot&\cdot&\cdot&\cdot&\cdot&\cdot&\cdot&\cdot&\cdot&\cdot&1&\cdot&\cdot&1&\cdot&\cdot&\cdot&\cdot&1\\
&4096_z'&\cdot&\cdot&\cdot&\cdot&\cdot&\cdot&\cdot&\cdot&\cdot&\cdot&\cdot&\cdot&1&1&\cdot&1&1&\cdot&1&\cdot\\
&2268_x'&\cdot&\cdot&\cdot&\cdot&\cdot&\cdot&\cdot&\cdot&\cdot&\cdot&\cdot&\cdot&\cdot&1&\cdot&\cdot&1&\cdot&\cdot&\cdot\\
&972_x'&\cdot&\cdot&\cdot&\cdot&\cdot&\cdot&\cdot&\cdot&\cdot&\cdot&\cdot&\cdot&\cdot&\cdot&1&1&\cdot&1&\cdot&1\\
&1344_x'&\cdot&\cdot&\cdot&\cdot&\cdot&\cdot&\cdot&\cdot&\cdot&\cdot&\cdot&\cdot&\cdot&\cdot&\cdot&1&1&1&\cdot&\cdot\\
&567_x'&\cdot&\cdot&\cdot&\cdot&\cdot&\cdot&\cdot&\cdot&\cdot&\cdot&\cdot&\cdot&\cdot&\cdot&\cdot&\cdot&1&\cdot&1&\cdot\\
&84_x'&\cdot&\cdot&\cdot&\cdot&\cdot&\cdot&\cdot&\cdot&\cdot&\cdot&\cdot&\cdot&\cdot&\cdot&\cdot&\cdot&\cdot&1&\cdot&1\\
&28_x'&\cdot&\cdot&\cdot&\cdot&\cdot&\cdot&\cdot&\cdot&\cdot&\cdot&\cdot&\cdot&\cdot&\cdot&\cdot&\cdot&\cdot&\cdot&1&\cdot\\
&1_x'&\cdot&\cdot&\cdot&\cdot&\cdot&\cdot&\cdot&\cdot&\cdot&\cdot&\cdot&\cdot&\cdot&\cdot&\cdot&\cdot&\cdot&\cdot&\cdot&1\\
\end{block}
\end{blockarray}
\]

\[
\begin{blockarray}{cccccccccccccccccccccc}
\begin{block}{cc(cccccccccccccccccccc)}
\star&1_x&1&-1&\cdot&\cdot&1&\cdot&\cdot&-1&1&\cdot&1&\cdot&-1&\cdot&\cdot&1&\cdot&-1&\cdot&1\\
\mathtt{(4)}&84_x&\cdot&1&\cdot&\cdot&-1&\cdot&\cdot&1&-1&\cdot&-1&-1&1&\cdot&1&-2&1&1&-1&-1\\
\star&28_x&\cdot&\cdot&1&-1&1&1&-1&-1&\cdot&1&1&1&-1&1&-1&1&-1&\cdot&1&\cdot\\
\mathtt{(4)}&567_x&\cdot&\cdot&\cdot&1&-1&-1&1&1&\cdot&-1&-2&-1&2&-2&1&-2&2&1&-2&-1\\
\mathtt{(4)}&1344_x&\cdot&\cdot&\cdot&\cdot&1&\cdot&-1&-1&1&1&1&1&-2&1&-1&2&-1&-1&2&1\\
\mathtt{(4)}&2268_x&\cdot&\cdot&\cdot&\cdot&\cdot&1&\cdot&-1&\cdot&1&1&\cdot&-1&1&-1&1&-1&\cdot&1&1\\
&972_x&\cdot&\cdot&\cdot&\cdot&\cdot&\cdot&1&\cdot&\cdot&-1&\cdot&-1&1&-1&1&-1&1&\cdot&-2&\cdot\\
&4096_z&\cdot&\cdot&\cdot&\cdot&\cdot&\cdot&\cdot&1&-1&-1&-1&\cdot&2&-1&1&-2&1&1&-2&-2\\
&4536_y&\cdot&\cdot&\cdot&\cdot&\cdot&\cdot&\cdot&\cdot&1&\cdot&\cdot&\cdot&-1&\cdot&\cdot&1&\cdot&-1&1&1\\
&2688_y&\cdot&\cdot&\cdot&\cdot&\cdot&\cdot&\cdot&\cdot&\cdot&1&\cdot&\cdot&-1&1&-1&1&-1&\cdot&2&1\\
&1134_y&\cdot&\cdot&\cdot&\cdot&\cdot&\cdot&\cdot&\cdot&\cdot&\cdot&1&\cdot&-1&1&\cdot&1&-1&-1&1&1\\
&168_y&\cdot&\cdot&\cdot&\cdot&\cdot&\cdot&\cdot&\cdot&\cdot&\cdot&\cdot&1&\cdot&\cdot&-1&1&-1&\cdot&1&\cdot\\
&4096_z'&\cdot&\cdot&\cdot&\cdot&\cdot&\cdot&\cdot&\cdot&\cdot&\cdot&\cdot&\cdot&1&-1&\cdot&-1&1&1&-2&-1\\
&2268_x'&\cdot&\cdot&\cdot&\cdot&\cdot&\cdot&\cdot&\cdot&\cdot&\cdot&\cdot&\cdot&\cdot&1&\cdot&\cdot&-1&\cdot&1&\cdot\\
&972_x'&\cdot&\cdot&\cdot&\cdot&\cdot&\cdot&\cdot&\cdot&\cdot&\cdot&\cdot&\cdot&\cdot&\cdot&1&-1&1&\cdot&-1&-1\\
&1344_x'&\cdot&\cdot&\cdot&\cdot&\cdot&\cdot&\cdot&\cdot&\cdot&\cdot&\cdot&\cdot&\cdot&\cdot&\cdot&1&-1&-1&1&1\\
&567_x'&\cdot&\cdot&\cdot&\cdot&\cdot&\cdot&\cdot&\cdot&\cdot&\cdot&\cdot&\cdot&\cdot&\cdot&\cdot&\cdot&1&\cdot&-1&\cdot\\
&84_x'&\cdot&\cdot&\cdot&\cdot&\cdot&\cdot&\cdot&\cdot&\cdot&\cdot&\cdot&\cdot&\cdot&\cdot&\cdot&\cdot&\cdot&1&\cdot&-1\\
&28_x'&\cdot&\cdot&\cdot&\cdot&\cdot&\cdot&\cdot&\cdot&\cdot&\cdot&\cdot&\cdot&\cdot&\cdot&\cdot&\cdot&\cdot&\cdot&1&\cdot\\
&1_x'&\cdot&\cdot&\cdot&\cdot&\cdot&\cdot&\cdot&\cdot&\cdot&\cdot&\cdot&\cdot&\cdot&\cdot&\cdot&\cdot&\cdot&\cdot&\cdot&1\\
\end{block}
\end{blockarray}
\]

\small
\begin{center}
\begin{picture}(490,100)
\put(0,80){\line(1,0){490}}
\put(0,80){\circle*{5}}
\put(65,80){\circle*{5}}
\put(115,80){\circle*{5}}
\put(145,80){\circle*{5}}
\put(175,80){\circle*{5}}
\put(195,80){\circle*{5}}
\put(245,80){\circle*{5}}
\put(295,80){\circle*{5}}
\put(315,80){\circle*{5}}
\put(345,80){\circle*{5}}
\put(375,80){\circle*{5}}
\put(425,80){\circle*{5}}
\put(490,80){\circle*{5}}
\put(-11,90){$-20$}
\put(55,90){$-8$}
\put(105,90){$-4$}
\put(135,90){$-2$}
\put(172,90){$0$}
\put(193,90){$1$}
\put(242,90){$4$}
\put(293,90){$7$}
\put(312,90){$8$}
\put(340,90){$10$}
\put(370,90){$12$}
\put(420,90){$16$}
\put(485,90){$28$}
\put(-2,65){$1_x$}
\put(60,65){$84_x$}
\put(60,50){$28_x$}
\put(105,65){$567_x$}
\put(132,65){$1344_x$}
\put(165,65){$972_x$}
\put(165,50){$2268_x$}
\put(190,65){$4096_z$}
\put(237,65){$168_y$}
\put(235,50){$1134_y$}
\put(235,35){$2688_y$}
\put(235,20){$4536_y$}
\put(280,65){$4096_z'$}
\put(310,65){$972_x'$}
\put(307,50){$2268_x'$}
\put(337,65){$1344_x'$}
\put(369,65){$567_x'$}
\put(419,65){$84_x'$}
\put(419,50){$28_x'$}
\put(485,65){$1_x'$}
\end{picture}
\end{center}
\normalsize

\textbf{\large Finite-dimensional modules, their characters and dimensions:}\\
\begin{align*}
\chi_{\el(1_x)}(t)&=t^{-20}+t^{20} + 8(t^{-19}+t^{19}) + 36(t^{-18}+t^{18}) + 120(t^{-17}+t^{17}) + 330(t^{-16}+t^{16}) + 792(t^{-15}+t^{15}) \\&\qquad+ 1716(t^{-14}+t^{14}) + 3432(t^{-13}+t^{13}) + 6435(t^{-12}+t^{12}) + 11440(t^{-11}+t^{11}) + 19448(t^{-10}+t^{10}) \\&\qquad + 31824(t^{-9}+t^9) + 50304(t^{-8}+t^8)+ 76848(t^{-7}+t^7) + 113256(t^{-6}+t^6) + 160464(t^{-5}+t^5) \\&\qquad +  217437(t^{-4}+t^4)+ 279576(t^{-3}+t^3) + 337900(t^{-2}+t^2) + 380264(t^{-1}+t) + 395874\\
\dim\el(1_x)&=3,779,136\\
\\
\chi_{\el(28_x)}(t)&=28(t^{-8}+t^8) + 224(t^{-7}+t^7) + 1008(t^{-6}+t^6) + 3360(t^{-5}+t^5) + 8673(t^{-4}+t^4) \\&\qquad+ 17640(t^{-3}+t^3) + 28980(t^{-2}+t^2) + 38808(t^{-1}+t) + 42750\\
\dim\el(28_x)&=240,192\\
\\
\end{align*}

\subsubsection{The block containing $\el(8_z)$}

\small

\[
\begin{blockarray}{cccccccccccccccccccccc}
\begin{block}{cc(cccccccccccccccccccc)}
\star&8_z&1&1&\cdot&\cdot&\cdot&\cdot&\cdot&\cdot&1&\cdot&\cdot&\cdot&\cdot&\cdot&\cdot&\cdot&\cdot&\cdot&\cdot&\cdot\\
\mathtt{(4)}&112_z&\cdot&1&1&1&\cdot&\cdot&\cdot&1&1&\cdot&\cdot&\cdot&\cdot&\cdot&\cdot&\cdot&\cdot&\cdot&\cdot&\cdot\\
\mathtt{(4)}&1008_z&\cdot&\cdot&1&\cdot&\cdot&1&1&1&\cdot&\cdot&\cdot&\cdot&\cdot&\cdot&\cdot&\cdot&\cdot&\cdot&\cdot&\cdot\\
\mathtt{(4)}&448_z&\cdot&\cdot&\cdot&1&\cdot&\cdot&\cdot&1&\cdot&1&\cdot&\cdot&\cdot&\cdot&\cdot&\cdot&\cdot&\cdot&\cdot&\cdot\\
\star&56_z&\cdot&\cdot&\cdot&\cdot&1&1&\cdot&\cdot&\cdot&\cdot&\cdot&1&\cdot&\cdot&\cdot&\cdot&\cdot&\cdot&\cdot&\cdot\\
\mathtt{(4)}&1296_z&\cdot&\cdot&\cdot&\cdot&\cdot&1&1&\cdot&\cdot&\cdot&\cdot&1&\cdot&\cdot&1&\cdot&\cdot&\cdot&\cdot&\cdot\\
&4096_x&\cdot&\cdot&\cdot&\cdot&\cdot&\cdot&1&1&\cdot&\cdot&1&\cdot&\cdot&\cdot&1&\cdot&\cdot&\cdot&\cdot&\cdot\\
&4536_z&\cdot&\cdot&\cdot&\cdot&\cdot&\cdot&\cdot&1&1&1&1&\cdot&1&\cdot&\cdot&\cdot&\cdot&\cdot&\cdot&\cdot\\
&1344_w&\cdot&\cdot&\cdot&\cdot&\cdot&\cdot&\cdot&\cdot&1&\cdot&\cdot&\cdot&1&\cdot&\cdot&\cdot&\cdot&\cdot&\cdot&1\\
&2016_w&\cdot&\cdot&\cdot&\cdot&\cdot&\cdot&\cdot&\cdot&\cdot&1&\cdot&\cdot&1&\cdot&\cdot&\cdot&1&\cdot&\cdot&\cdot\\
&7168_w&\cdot&\cdot&\cdot&\cdot&\cdot&\cdot&\cdot&\cdot&\cdot&\cdot&1&\cdot&1&1&1&\cdot&\cdot&\cdot&\cdot&\cdot\\
&448_w&\cdot&\cdot&\cdot&\cdot&\cdot&\cdot&\cdot&\cdot&\cdot&\cdot&\cdot&1&\cdot&\cdot&1&\cdot&\cdot&1&\cdot&\cdot\\
&4536_z'&\cdot&\cdot&\cdot&\cdot&\cdot&\cdot&\cdot&\cdot&\cdot&\cdot&\cdot&\cdot&1&1&\cdot&\cdot&1&\cdot&1&1\\
&4096_x'&\cdot&\cdot&\cdot&\cdot&\cdot&\cdot&\cdot&\cdot&\cdot&\cdot&\cdot&\cdot&\cdot&1&1&1&\cdot&\cdot&1&\cdot\\
&1296_z'&\cdot&\cdot&\cdot&\cdot&\cdot&\cdot&\cdot&\cdot&\cdot&\cdot&\cdot&\cdot&\cdot&\cdot&1&1&\cdot&1&\cdot&\cdot\\
&1008_z'&\cdot&\cdot&\cdot&\cdot&\cdot&\cdot&\cdot&\cdot&\cdot&\cdot&\cdot&\cdot&\cdot&\cdot&\cdot&1&\cdot&\cdot&1&\cdot\\
&448_z'&\cdot&\cdot&\cdot&\cdot&\cdot&\cdot&\cdot&\cdot&\cdot&\cdot&\cdot&\cdot&\cdot&\cdot&\cdot&\cdot&1&\cdot&1&\cdot\\
&56_z'&\cdot&\cdot&\cdot&\cdot&\cdot&\cdot&\cdot&\cdot&\cdot&\cdot&\cdot&\cdot&\cdot&\cdot&\cdot&\cdot&\cdot&1&\cdot&\cdot\\
&112_z'&\cdot&\cdot&\cdot&\cdot&\cdot&\cdot&\cdot&\cdot&\cdot&\cdot&\cdot&\cdot&\cdot&\cdot&\cdot&\cdot&\cdot&\cdot&1&1\\
&8_z'&\cdot&\cdot&\cdot&\cdot&\cdot&\cdot&\cdot&\cdot&\cdot&\cdot&\cdot&\cdot&\cdot&\cdot&\cdot&\cdot&\cdot&\cdot&\cdot&1\\
\end{block}
\end{blockarray}
\]

\[
\begin{blockarray}{cccccccccccccccccccccc}
\begin{block}{cc(cccccccccccccccccccc)}
\star&8_z&1&-1&1&1&\cdot&-1&\cdot&-1&1&\cdot&1&1&-1&\cdot&-1&1&1&\cdot&-1&1\\
\mathtt{(4)}&112_z&\cdot&1&-1&-1&\cdot&1&\cdot&1&-2&\cdot&-1&-1&2&-1&2&-1&-2&-1&2&-2\\
\mathtt{(4)}&1008_z&\cdot&\cdot&1&\cdot&\cdot&-1&\cdot&-1&1&1&1&1&-2&1&-2&1&1&1&-1&2\\
\mathtt{(4)}&448_z&\cdot&\cdot&\cdot&1&\cdot&\cdot&\cdot&-1&1&\cdot&1&\cdot&-1&\cdot&-1&1&1&1&-1&1\\
\star&56_z&\cdot&\cdot&\cdot&\cdot&1&-1&1&-1&1&1&\cdot&\cdot&-1&1&-1&\cdot&\cdot&1&\cdot&\cdot\\
\mathtt{(4)}&1296_z&\cdot&\cdot&\cdot&\cdot&\cdot&1&-1&1&-1&-1&\cdot&-1&1&-1&2&-1&\cdot&-1&1&-1\\
&4096_x&\cdot&\cdot&\cdot&\cdot&\cdot&\cdot&1&-1&1&1&\cdot&\cdot&-1&1&-2&1&\cdot&2&-1&1\\
&4536_z&\cdot&\cdot&\cdot&\cdot&\cdot&\cdot&\cdot&1&-1&-1&-1&\cdot&2&-1&2&-1&-1&-2&1&-2\\
&1344_w&\cdot&\cdot&\cdot&\cdot&\cdot&\cdot&\cdot&\cdot&1&\cdot&\cdot&\cdot&-1&1&-1&\cdot&1&1&-1&1\\
&2016_w&\cdot&\cdot&\cdot&\cdot&\cdot&\cdot&\cdot&\cdot&\cdot&1&\cdot&\cdot&-1&1&-1&\cdot&\cdot&1&\cdot&1\\
&7168_w&\cdot&\cdot&\cdot&\cdot&\cdot&\cdot&\cdot&\cdot&\cdot&\cdot&1&\cdot&-1&\cdot&-1&1&1&1&-1&2\\
&448_w&\cdot&\cdot&\cdot&\cdot&\cdot&\cdot&\cdot&\cdot&\cdot&\cdot&\cdot&1&\cdot&\cdot&-1&1&\cdot&\cdot&-1&1\\
&4536_z'&\cdot&\cdot&\cdot&\cdot&\cdot&\cdot&\cdot&\cdot&\cdot&\cdot&\cdot&\cdot&1&-1&1&\cdot&-1&-1&1&-2\\
&4096_x'&\cdot&\cdot&\cdot&\cdot&\cdot&\cdot&\cdot&\cdot&\cdot&\cdot&\cdot&\cdot&\cdot&1&-1&\cdot&\cdot&1&-1&1\\
&1296_z'&\cdot&\cdot&\cdot&\cdot&\cdot&\cdot&\cdot&\cdot&\cdot&\cdot&\cdot&\cdot&\cdot&\cdot&1&-1&\cdot&-1&1&-1\\
&1008_z'&\cdot&\cdot&\cdot&\cdot&\cdot&\cdot&\cdot&\cdot&\cdot&\cdot&\cdot&\cdot&\cdot&\cdot&\cdot&1&\cdot&\cdot&-1&1\\
&448_z'&\cdot&\cdot&\cdot&\cdot&\cdot&\cdot&\cdot&\cdot&\cdot&\cdot&\cdot&\cdot&\cdot&\cdot&\cdot&\cdot&1&\cdot&-1&1\\
&56_z'&\cdot&\cdot&\cdot&\cdot&\cdot&\cdot&\cdot&\cdot&\cdot&\cdot&\cdot&\cdot&\cdot&\cdot&\cdot&\cdot&\cdot&1&\cdot&\cdot\\
&112_z'&\cdot&\cdot&\cdot&\cdot&\cdot&\cdot&\cdot&\cdot&\cdot&\cdot&\cdot&\cdot&\cdot&\cdot&\cdot&\cdot&\cdot&\cdot&1&-1\\
&8_z'&\cdot&\cdot&\cdot&\cdot&\cdot&\cdot&\cdot&\cdot&\cdot&\cdot&\cdot&\cdot&\cdot&\cdot&\cdot&\cdot&\cdot&\cdot&\cdot&1\\
\end{block}
\end{blockarray}
\]

\small

\begin{center}
\begin{picture}(480,100)
\put(0,80){\line(1,0){480}}
\put(0,80){\circle*{5}}
\put(60,80){\circle*{5}}
\put(120,80){\circle*{5}}
\put(160,80){\circle*{5}}
\put(180,80){\circle*{5}}
\put(200,80){\circle*{5}}
\put(240,80){\circle*{5}}
\put(280,80){\circle*{5}}
\put(300,80){\circle*{5}}
\put(320,80){\circle*{5}}
\put(360,80){\circle*{5}}
\put(420,80){\circle*{5}}
\put(480,80){\circle*{5}}
\put(-11,90){$-14$}
\put(52,90){$-8$}
\put(112,90){$-2$}
\put(157,90){$0$}
\put(177,90){$1$}
\put(197,90){$2$}
\put(237,90){$4$}
\put(277,90){$6$}
\put(297,90){$7$}
\put(317,90){$8$}
\put(355,90){$10$}
\put(415,90){$16$}
\put(475,90){$22$}
\put(-2,65){$8_z$}
\put(53,65){$112_z$}
\put(110,65){$1008_z$}
\put(112,50){$448_z$}
\put(113,35){$56_z$}
\put(143,65){$1296_z$}
\put(170,65){$4096_x$}
\put(196,65){$4536_z$}
\put(230,65){$1344_w$}
\put(230,50){$2016_w$}
\put(230,35){$7168_w$}
\put(230,20){$448_w$}
\put(263,65){$4536_z'$}
\put(290,65){$4096_x'$}
\put(317,65){$1296_z'$}
\put(350,65){$1008_z'$}
\put(352,50){$448_z'$}
\put(353,35){$56_z'$}
\put(413,65){$112_z'$}
\put(477,65){$8_z'$}
\end{picture}
\end{center}

\normalsize
\textbf{\large Finite-dimensional modules, their characters and dimensions:}\\
\begin{align*}
\chi_{\el(8_z)}(t)&=8(t^{-14}+t^{14}) + 64(t^{-13}+t^{13}) + 288(t^{-12}+t^{12}) + 960(t^{-11}+t^{11}) + 2640(t^{-10}+t^{10})\\&\qquad + 6336(t^{-9}+t^9) + 13616(t^{-8}+t^8) + 26560(t^{-7}+t^7) + 47448(t^{-6}+t^6) + 78080(t^{-5}+t^5)\\&\qquad + 118624(t^{-4}+t^4) + 165888(t^{-3}+t^3) + 212368(t^{-2}+t^2) + 247424(t^{-1}+t) + 260640\\
\dim\el(8_z)&=2,101,248\\
\\
\chi_{\el(56_z)}(t)&=56(t^{-2}+t^2) + 448(t^{-1}+t) + 720\\
\dim\el(56_z)&=1728\\
\end{align*}

\subsection{{\Large$\mbf{c=\frac{1}{4}}$}}
The principal block is large and not studied in this paper. There are also four smaller blocks, of which two are self-dual and two are dual to each other under tensoring the characters labeling irreducibles with the sign representation. The decomposition matrices of these four blocks are worked out below. \\

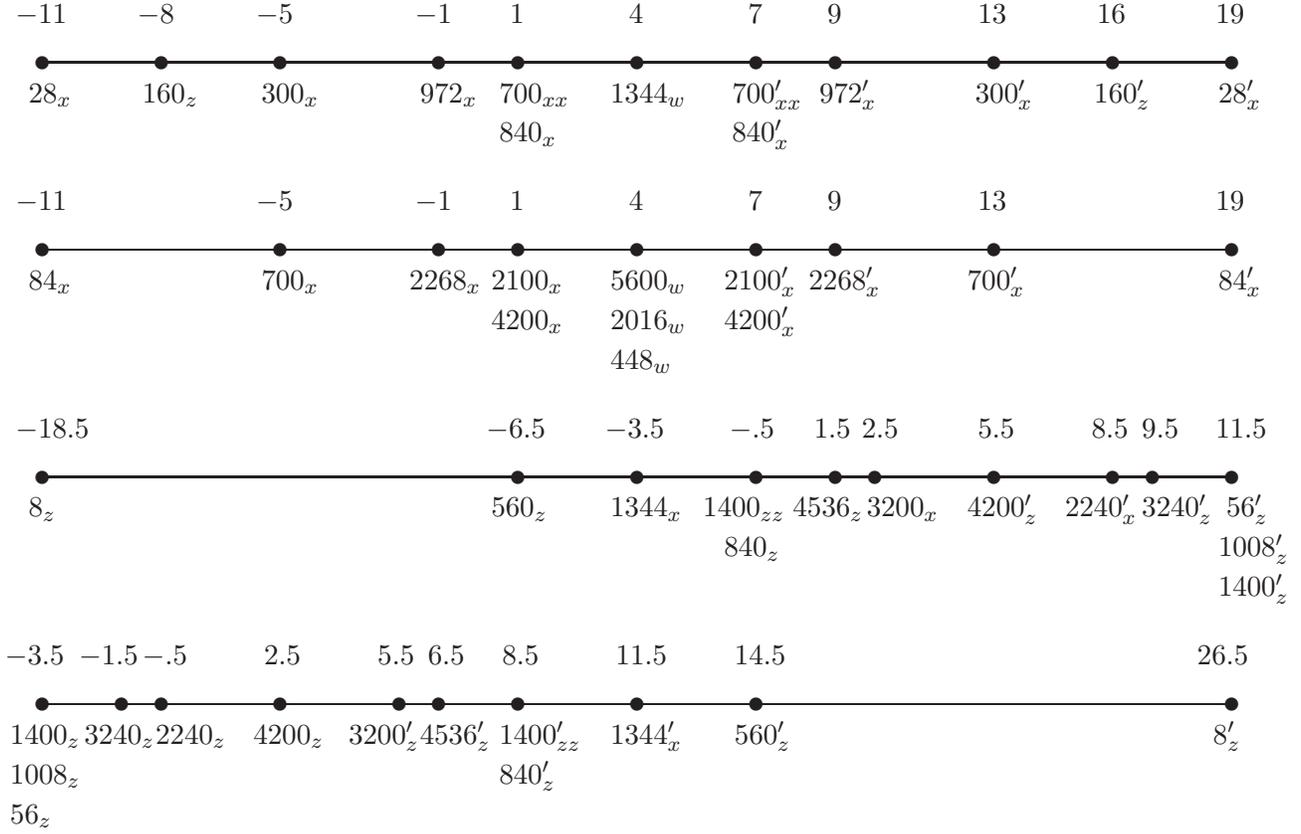
\begin{center}
\begin{picture}(450,70)
\put(0,50){\line(1,0){450}}
\put(0,50){\circle*{5}}
\put(45,50){\circle*{5}}
\put(90,50){\circle*{5}}
\put(150,50){\circle*{5}}
\put(180,50){\circle*{5}}
\put(225,50){\circle*{5}}
\put(270,50){\circle*{5}}
\put(300,50){\circle*{5}}
\put(360,50){\circle*{5}}
\put(405,50){\circle*{5}}
\put(450,50){\circle*{5}}
\put(-10,65){$-11$}
\put(36,65){$-8$}
\put(81,65){$-5$}
\put(141,65){$-1$}
\put(177,65){$1$}
\put(222,65){$4$}
\put(267,65){$7$}
\put(297,65){$9$}
\put(354,65){$13$}
\put(399,65){$16$}
\put(444,65){$19$}
\put(-5,35){$28_x$}
\put(38,35){$160_z$}
\put(83,35){$300_x$}
\put(143,35){$972_x$}
\put(173,35){$700_{xx}$}
\put(173,20){$840_x$}
\put(215,35){$1344_w$}
\put(261,35){$700_{xx}'$}
\put(261,20){$840_x'$}
\put(294,35){$972_x'$}
\put(353,35){$300_x'$}
\put(398,35){$160_z'$}
\put(445,35){$28_x'$}
\end{picture}

\begin{picture}(450,85)
\put(0,65){\line(1,0){450}}
\put(0,65){\circle*{5}}
\put(90,65){\circle*{5}}
\put(150,65){\circle*{5}}
\put(180,65){\circle*{5}}
\put(225,65){\circle*{5}}
\put(270,65){\circle*{5}}
\put(300,65){\circle*{5}}
\put(360,65){\circle*{5}}
\put(450,65){\circle*{5}}
\put(-10,80){$-11$}
\put(81,80){$-5$}
\put(141,80){$-1$}
\put(177,80){$1$}
\put(222,80){$4$}
\put(267,80){$7$}
\put(297,80){$9$}
\put(354,80){$13$}
\put(444,80){$19$}
\put(-5,50){$84_x$}
\put(83,50){$700_x$}
\put(139,50){$2268_x$}
\put(170,50){$2100_x$}
\put(170,35){$4200_x$}
\put(215,50){$5600_w$}
\put(215,35){$2016_w$}
\put(215,20){$448_w$}
\put(258,50){$2100_x'$}
\put(258,35){$4200_x'$}
\put(290,50){$2268_x'$}
\put(350,50){$700_x'$}
\put(445,50){$84_x'$}
\end{picture}

\begin{picture}(450,85)
\put(0,65){\line(1,0){450}}
\put(0,65){\circle*{5}}
\put(180,65){\circle*{5}}
\put(225,65){\circle*{5}}
\put(270,65){\circle*{5}}
\put(300,65){\circle*{5}}
\put(315,65){\circle*{5}}
\put(360,65){\circle*{5}}
\put(405,65){\circle*{5}}
\put(420,65){\circle*{5}}
\put(450,65){\circle*{5}}
\put(-10,80){$-18.5$}
\put(168,80){$-6.5$}
\put(213,80){$-3.5$}
\put(260,80){$-.5$}
\put(292,80){$1.5$}
\put(310,80){$2.5$}
\put(354,80){$5.5$}
\put(397,80){$8.5$}
\put(416,80){$9.5$}
\put(444,80){$11.5$}
\put(-5,50){$8_z$}
\put(170,50){$560_z$}
\put(215,50){$1344_x$}
\put(250,50){$1400_{zz}$}
\put(258,35){$840_z$}
\put(284,50){$4536_z$}
\put(312,50){$3200_x$}
\put(350,50){$4200_z'$}
\put(387,50){$2240_x'$}
\put(416,50){$3240_z'$}
\put(448,50){$56_z'$}
\put(445,35){$1008_z'$}
\put(445,20){$1400_z'$}
\end{picture}

\begin{picture}(450,85)
\put(0,65){\line(1,0){450}}
\put(0,65){\circle*{5}}
\put(30,65){\circle*{5}}
\put(45,65){\circle*{5}}
\put(90,65){\circle*{5}}
\put(135,65){\circle*{5}}
\put(150,65){\circle*{5}}
\put(180,65){\circle*{5}}
\put(225,65){\circle*{5}}
\put(270,65){\circle*{5}}
\put(450,65){\circle*{5}}
\put(-14,80){$-3.5$}
\put(14,80){$-1.5$}
\put(38,80){$-.5$}
\put(84,80){$2.5$}
\put(127,80){$5.5$}
\put(146,80){$6.5$}
\put(174,80){$8.5$}
\put(217,80){$11.5$}
\put(262,80){$14.5$}
\put(437,80){$26.5$}
\put(-12,50){$1400_z$}
\put(-12,35){$1008_z$}
\put(-12,20){$56_z$}
\put(16,50){$3240_z$}
\put(43,50){$2240_z$}
\put(80,50){$4200_z$}
\put(116,50){$3200_z'$}
\put(143,50){$4536_z'$}
\put(173,50){$1400_{zz}'$}
\put(173,35){$840_z'$}
\put(215,50){$1344_x'$}
\put(262,50){$560_z'$}
\put(443,50){$8_z'$}
\end{picture}
\end{center}

The block containing $\el(28_x)$ is self-dual and contains $3$ irreps with dimension $4$ support and $2$ irreps with dimension $5$ support. The irreps with dimension $4$ support have decompositions invariant under tensoring the lowest weights of the Vermas involved with the sign rep, and consequently their characters 
$$\chi_{\el(\tau)}(t)=t^{h_c(\tau)}\frac{P_\tau(t)}{(1-t)^4}$$
where $P_\tau(t)$ is a polynomial of degree $N:=2(4-h_c(\tau))$ are invariant under replacing $t^k$ with $t^{N-k}$ for each $k=0,...,N$, i.e. $$t^{\frac{-N}{2}}P_\tau(t)=t^{\frac{N}{2}}P_\tau(-t)$$
All but the first five columns appear in the Hecke algebra decomposition matrix for this block, and the entries for the first five columns follow immediately by applying Lemmas \ref{RR} and \ref{dim Hom}.

\[
\begin{blockarray}{ccccccccccccccc}
&&28_x&160_z&300_x&972_x&700_{xx}&840_x&1344_w&700_{xx}'&840_x'&972_x'&300_x'&160_z'&28_x'\\
\begin{block}{cc(ccccccccccccc)}
\mathtt{(4)}&28_x&1&\cdot&\cdot&\cdot&1&\cdot&\cdot&\cdot&1&\cdot&\cdot&\cdot&\cdot\\
\mathtt{(4)}&160_z&\cdot&1&\cdot&1&1&\cdot&1&\cdot&\cdot&\cdot&\cdot&\cdot&\cdot\\
\mathtt{(4)}&300_x&\cdot&\cdot&1&1&\cdot&1&\cdot&\cdot&\cdot&\cdot&\cdot&\cdot&\cdot\\
\mathtt{(5)}&972_x&\cdot&\cdot&\cdot&1&\cdot&1&1&1&\cdot&\cdot&\cdot&\cdot&\cdot\\
\mathtt{(5)}&700_{xx}&\cdot&\cdot&\cdot&\cdot&1&\cdot&1&\cdot&1&1&\cdot&\cdot&\cdot\\
&840_x&\cdot&\cdot&\cdot&\cdot&\cdot&1&\cdot&1&\cdot&\cdot&\cdot&\cdot&1\\
&1344_w&\cdot&\cdot&\cdot&\cdot&\cdot&\cdot&1&1&\cdot&1&\cdot&1&\cdot\\
&700_{xx}'&\cdot&\cdot&\cdot&\cdot&\cdot&\cdot&\cdot&1&\cdot&\cdot&\cdot&1&1\\
&840_x'&\cdot&\cdot&\cdot&\cdot&\cdot&\cdot&\cdot&\cdot&1&1&1&\cdot&\cdot\\
&972_x'&\cdot&\cdot&\cdot&\cdot&\cdot&\cdot&\cdot&\cdot&\cdot&1&1&1&\cdot\\
&300_x'&\cdot&\cdot&\cdot&\cdot&\cdot&\cdot&\cdot&\cdot&\cdot&\cdot&1&\cdot&\cdot\\
&160_z'&\cdot&\cdot&\cdot&\cdot&\cdot&\cdot&\cdot&\cdot&\cdot&\cdot&\cdot&1&\cdot\\
&28_x'&\cdot&\cdot&\cdot&\cdot&\cdot&\cdot&\cdot&\cdot&\cdot&\cdot&\cdot&\cdot&1\\
\end{block}\\
\begin{block}{cc(ccccccccccccc)}
\mathtt{(4)}&28_x&1&\cdot&\cdot&\cdot&-1&\cdot&1&-1&\cdot&\cdot&\cdot&\cdot&1\\
\mathtt{(4)}&160_z&\cdot&1&\cdot&-1&-1&1&1&-1&1&-1&\cdot&1&\cdot\\
\mathtt{(4)}&300_x&\cdot&\cdot&1&-1&\cdot&\cdot&1&\cdot&\cdot&-1&1&\cdot&\cdot\\
\mathtt{(5)}&972_x&\cdot&\cdot&\cdot&1&\cdot&-1&-1&1&\cdot&1&-1&-1&\cdot\\
\mathtt{(5)}&700_{xx}&\cdot&\cdot&\cdot&\cdot&1&\cdot&-1&1&-1&1&\cdot&-1&-1\\
&840_x&\cdot&\cdot&\cdot&\cdot&\cdot&1&\cdot&-1&\cdot&\cdot&\cdot&1&\cdot\\
&1344_w&\cdot&\cdot&\cdot&\cdot&\cdot&\cdot&1&-1&\cdot&-1&1&1&1\\
&700_{xx}'&\cdot&\cdot&\cdot&\cdot&\cdot&\cdot&\cdot&1&\cdot&\cdot&\cdot&-1&-1\\
&840_x'&\cdot&\cdot&\cdot&\cdot&\cdot&\cdot&\cdot&\cdot&1&-1&\cdot&1&\cdot\\
&972_x'&\cdot&\cdot&\cdot&\cdot&\cdot&\cdot&\cdot&\cdot&\cdot&1&-1&-1&\cdot\\
&300_x'&\cdot&\cdot&\cdot&\cdot&\cdot&\cdot&\cdot&\cdot&\cdot&\cdot&1&\cdot&\cdot\\
&160_z'&\cdot&\cdot&\cdot&\cdot&\cdot&\cdot&\cdot&\cdot&\cdot&\cdot&\cdot&1&\cdot\\
&28_x'&\cdot&\cdot&\cdot&\cdot&\cdot&\cdot&\cdot&\cdot&\cdot&\cdot&\cdot&\cdot&1\\
\end{block}
\end{blockarray}
\]

The block containing $84_x$ is also self-dual.

\[
\begin{blockarray}{ccccccccccccccc}
&&84_x&700_x&2268_x&2100_x&4200_x&5600_w&2016_w&448_w&2100_x'&4200_x'&2268_x'&700_x'&84_x'\\
\begin{block}{cc(ccccccccccccc)}
\mathtt{(4)}&84_x&1&1&\cdot&\cdot&\cdot&\cdot&1&\cdot&\cdot&\cdot&\cdot&\cdot&\cdot\\
\mathtt{(5)}&700_x&\cdot&1&1&\cdot&1&\cdot&1&\cdot&\cdot&\cdot&\cdot&\cdot&\cdot\\
\mathtt{(4)}&2268_x&\cdot&\cdot&1&1&1&1&\cdot&\cdot&\cdot&\cdot&\cdot&\cdot&\cdot\\
\mathtt{(5)}&2100_x&\cdot&\cdot&\cdot&1&\cdot&1&\cdot&1&1&\cdot&\cdot&\cdot&\cdot\\
&4200_x&\cdot&\cdot&\cdot&\cdot&1&1&1&\cdot&\cdot&1&\cdot&\cdot&\cdot\\
&5600_w&\cdot&\cdot&\cdot&\cdot&\cdot&1&\cdot&\cdot&1&1&1&\cdot&\cdot\\
&2016_w&\cdot&\cdot&\cdot&\cdot&\cdot&\cdot&1&\cdot&\cdot&1&\cdot&\cdot&1\\
\mathtt{(4)}&448_w&\cdot&\cdot&\cdot&\cdot&\cdot&\cdot&\cdot&1&1&\cdot&\cdot&\cdot&\cdot\\
&2100_x'&\cdot&\cdot&\cdot&\cdot&\cdot&\cdot&\cdot&\cdot&1&\cdot&1&\cdot&\cdot\\
&4200_x'&\cdot&\cdot&\cdot&\cdot&\cdot&\cdot&\cdot&\cdot&\cdot&1&1&1&1\\
&2268_x'&\cdot&\cdot&\cdot&\cdot&\cdot&\cdot&\cdot&\cdot&\cdot&\cdot&1&1&\cdot\\
&700_x'&\cdot&\cdot&\cdot&\cdot&\cdot&\cdot&\cdot&\cdot&\cdot&\cdot&\cdot&1&1\\
&84_x'&\cdot&\cdot&\cdot&\cdot&\cdot&\cdot&\cdot&\cdot&\cdot&\cdot&\cdot&\cdot&1\\
\end{block}\\
\begin{block}{cc(ccccccccccccc)}
\mathtt{(4)}&84_x&1&-1&1&-1&\cdot&\cdot&\cdot&1&\cdot&\cdot&\cdot&\cdot&\cdot\\
\mathtt{(5)}&700_x& \cdot & 1& -1 & 1 & \cdot & \cdot& -1& -1&  \cdot&  1& -1&  \cdot&  \cdot\\
\mathtt{(4)}&2268_x& \cdot & \cdot & 1 &-1 &-1 & 1 & 1 & 1 &-1 &-1 & 1 & \cdot & \cdot\\
\mathtt{(5)}&2100_x& \cdot & \cdot & \cdot & 1 & \cdot &-1 & \cdot &-1 & 1 & 1 &-1 & \cdot& -1\\
&4200_x& \cdot & \cdot & \cdot & \cdot & 1 &-1& -1&  \cdot&  1&  1 &-1 & \cdot & \cdot\\
&5600_w& \cdot & \cdot & \cdot & \cdot & \cdot & 1 & \cdot & \cdot &-1 &-1 & 1 & \cdot & 1\\
&2016_w& \cdot & \cdot & \cdot & \cdot & \cdot & \cdot & 1 & \cdot & \cdot &-1 & 1 & \cdot & \cdot\\
\mathtt{(4)}&448_w& \cdot & \cdot & \cdot & \cdot&  \cdot & \cdot & \cdot & 1 &-1 & \cdot & 1 &-1 & 1\\
&2100_x'& \cdot & \cdot & \cdot & \cdot & \cdot&  \cdot & \cdot & \cdot &1 & \cdot &-1 & 1 &-1\\
&4200_x'& \cdot & \cdot & \cdot & \cdot & \cdot & \cdot & \cdot & \cdot & \cdot&  1 &-1 & \cdot &-1\\
&2268_x'& \cdot & \cdot & \cdot & \cdot & \cdot & \cdot & \cdot & \cdot & \cdot & \cdot & 1 &-1&  1\\
&700_x'& \cdot  &\cdot & \cdot & \cdot & \cdot & \cdot & \cdot & \cdot & \cdot & \cdot & \cdot & 1 &-1\\
&84_x'& \cdot & \cdot & \cdot & \cdot & \cdot & \cdot & \cdot  &\cdot  &\cdot  &\cdot  &\cdot  &\cdot & 1\\
\end{block}
\end{blockarray}
\]

Of the five columns killed by KZ, there are six entries of the decomposition matrix which cannot be recovered immediately by Lemmas \ref{RR} and \ref{dim Hom}. These are $[\m(84_x):\el(2268_x)],\;[\m(84_x):\el(2100_x)],\;[\m(84_x):\el(448_w)],\;[\m(700_x):\el(2100_x)],\;[\m(700_x):\el(448_w)],$ and $[\m(2268_x):\el(448_w)]$. We will use induction and restriction to the rational Cherednik algebra of $E_7$ at parameter $c=\frac{1}{4}$ to find the missing entries.

Let $\alpha=[\m(2268_x):\el(448_w)]$ and write:
\begin{align*}
\el(2268_x)&=\m(2268_x)-\m(2100_x)-\m(4200_x)+\m(5600_w)+\m(2016_w)+(1-\alpha)\m(448_w)\\
&\qquad+(\alpha-1)\m(2100_x')-\m(4200_x')+(1-\alpha)\m(2268_x')+\alpha\m(700_x')-\alpha\m(84_x')\\
\end{align*}
From the decompositions of $\el(448_w)=\m(448_w)-\m(2100_x')+\m(2268_x')...$ $\el(2100_x')=\m(2100_x')-\m(2268_x')...$ it follows that $\dim\Hom(\m(2268_x'),\m(448_w))=0$. Then by Lemma \ref{dim Hom}, $\dim\Hom(\m(448_w),\m(2268_x))=0$ and so $\alpha$ is either $0$ or $1$. If $\alpha=1$ then the following expression $M$ from the Grothendieck group would have to be the decomposition of $\el(448_w)$:

$$M:=\m(2268_x)-\m(2100_x)-\m(4200_x)+\m(5600_w)+\m(2016_w)-\m(4200_x')+\m(700_x')-\m(84_x')$$

Restricting $M$ to $H_{\frac{1}{4}}(E_7)$ we'd then have to get an expression for a module in the Category $\oh$.

\begin{align*}
\Res M&=\m(120_a)+\m(56_a')-\m(35_a')-\m(105_b)-\m(210_a)-\m(189_b')+\m(405_a)+\m(189_a')\\&\qquad+\m(70_a)+\m(70_a')-\m(210_a')-\m(315_a)+\m(120_a')+\m(56_a)-\m(21_b)-\m(1_a')\\
\end{align*}
The Vermas of minimal $h_c$-lowest-weight in this expression are $\m(120_a)$ and $\m(56_a')$ and these belong to different blocks, so $\el(120_a)$ and $\el(56_a')$ must both appear in the composition series of $\Res M$. Subtracting off these simples we have:
\begin{align*}
\Res M-\el(120_a)-\el(56_a')&=-\m(35_a')-\m(189_a)+\m(336_a')+\m(336_a)-\m(405_a')\\&\qquad-\m(315_a)+\m(105_b')+\m(189_b)-\m(21_b)-\m(1_a')\\
\end{align*}
But the Verma of minimal $h_c$-lowest-weight in this expression is $\m(35_a')$ and it appears with a minus sign, so $\Res M-\el(120_a)-\el(56_a')$ isn't the expression for a module in Category $\oh$ and thus $M$ couldn't have been a module either. Therefore $M\neq\el(2268_x)$, $\alpha\neq1$, and so $\alpha=0$.

The decomposition of $\el(700_x)$ is the next question. Either $[\m(700_x):\el(2100_x)]=0$ or $1$. If $[\m(700_x):\el(2100_x)]=1$ then restricting $\el(700_x)$ to the block containing $\el(1_a)$ in $\oh_{\frac{1}{4}}(E_7)$ produces the term $-\m(210_a)$ which neither belongs to any simple to its left in the restriction of this hypothetical $\el(700_x)$ nor is canceled by any terms to its right; this can't be a module, and so $[\m(700_x):\el(2100_x)]=0$.  The possible characters of $\el(700_x)$, depending on the value of $\alpha:=[\m(700_x):\el(448_w)]$, then all agree that the support of $\el(700_x)=5$; the decomposition of $\el(700_x)$ is:
\begin{align*}
\el(700_x)&=\m(700_x)-\m(2268_x)+\m(2100_x)-\m(2016_w)-(1+\alpha)\m(448_w)+\alpha\m(2100_x')\\
&\qquad+\m(4200_x')-(1+\alpha)\m(2268_x')+\alpha\m(700_x')-\alpha\m(84_x')\\
\end{align*}





 The restriction of $\el(700_x)$ to the block containing $\el(1_a)$ is:
\begin{align*}
\Res\el(700_x)|_{\B(1_a)}&=\m(105_b)-\m(405_a)+\m(336_a')-\m(70_a)-\m(35_a)+\m(189_b)-\m(120_a')\\
&\qquad-\alpha\left(\m(189_a)-\m(336_a')+\m(315_a)-\m(189_b)+\m(21_b)\right)\\
&=\el(105_b)-\alpha\el(189_a)\\
\end{align*}
Therefore $\alpha=[\m(700_x):\el(448_w)]=0$.

\newpage

The blocks containing $\el(35_b)$, $\el(35_b')$, and $\el(1_a)$ in $\oh_{\frac{1}{4}}(E_7)$ can be used to find the expression for $[\el(84_x)]$ in the basis of Vermas. Set $c_\tau:=[\el(84_x):\m(\tau)]$. Consider $\Ind\el(35_b)|_{\B(700_x)}$ precomposed with projection to the principal block of $\oh_{\frac{1}{4}}(E_7)$:
\begin{align*}
\Ind\el(35_b)|_{\B(700_x)}&=35_b\uparrow-210_b\uparrow+280_b'\uparrow-105a\uparrow\\
&=84_x+700_x-4200_x-2016_w+5600_w+4200_x'-2100_x'-2268_x'\\
\end{align*}
and
\begin{align*}
\Res(\Ind\el(35_b)|_{\B(700_x)})|_{\B(1_a)}&=1_a+56_a'-210_a-70_a+189_a+189_b-35_a-120_a'\\
&=\el(1_a)+\el(56_a')+2\el(105_b)\\
\end{align*}
Since $\el(84_x)=\m(84_x)-\m(700_x)+...$, $\Ind\el(35_b)|_{\B(700_x)}$ contains $\el(84_x)+2\el(700_x)$. The rep $56_a'$ appears only in $\Res700_x$ and $\Res2268_x$ among the $\tau\in\Irr E_8$ in the same block as $\el(84_x)$. Also, $[\el(1_a):\m(56_a')]=0$, and $\Res\el(700_x)|_{\B(1_a)}=\el(105_b)$, so it must be that $\Res\el(84_x)|_{\B(1_a)}=\el(1_a)$. We find immediately that $c_{2268_x}=1$. Next, the restriction of $\Ind\el(35_b)|_{\B(700_x)}$ to the block of $\oh_{\frac{1}{4}}(E_7)$ containing $\el(35_b')$ is $0$. Therefore $\Res\el(84_x)|_{\B(35_b')}=0$. It then follows from computing restrictions of the $\tau$ at hand that $c_{2100_x}=-c_{2268_x}$, so $c_{2100_x}=-1$. 

Thus $[\m(84_x):\el(2268_x)]=[\m(84_x):\el(2100_x]=0$, and the only missing entry of the decomposition matrix that's left to find is $[\m(84_x):\el(448_w)]$. Applying $\Res|_{\B(1_a)}$ to the decomposition of $\el(84_x)$ up through the $h_c$-weight-$4$ irreps except $\m(448_w)$, the resulting expression turns out to include the term $-\m(189_a)$. This term does not belong to $\el(1_a)$ so must be canceled out. Only $\m(448_w)$ remains free to do this as no $\tau$ to the right of $448_w$ has $189_a$ in its restriction, so $[\el(84_x):\m(448_w)]=1$ and $[\m(84_x):\el(448_w)]=0$. This completes the decomposition matrix of the block of $\el(84_x)$.

The block containing the simple module with lowest weight the standard representation has the following decomposition matrix:
\[
\begin{blockarray}{ccccccccccccccc}
\begin{block}{cc(ccccccccccccc)}
\mathtt{(4)}&8_z&1&\cdot&\cdot&1&\cdot&\cdot&\cdot&\cdot&1&\cdot&\cdot&\cdot&\cdot\\
\mathtt{(4)}&560_z&\cdot&1&1&1&\cdot&1&\cdot&\cdot&\cdot&\cdot&\cdot&\cdot&\cdot\\
\mathtt{(5)}&1344_x&\cdot&\cdot&1&\cdot&1&1&1&\cdot&\cdot&\cdot&\cdot&\cdot&\cdot\\
\mathtt{(5)}&1400_{zz}&\cdot&\cdot&\cdot&1&\cdot&1&\cdot&\cdot&1&1&\cdot&\cdot&\cdot\\
\mathtt{(4)}&840_z&\cdot&\cdot&\cdot&\cdot&1&\cdot&1&\cdot&\cdot&\cdot&\cdot&\cdot&1\\
&4536_z&\cdot&\cdot&\cdot&\cdot&\cdot&1&1&1&\cdot&1&\cdot&\cdot&\cdot\\
&3200_x&\cdot&\cdot&\cdot&\cdot&\cdot&\cdot&1&1&\cdot&\cdot&\cdot&\cdot&1\\
&4200_z'&\cdot&\cdot&\cdot&\cdot&\cdot&\cdot&\cdot&1&\cdot&1&\cdot&1&\cdot\\
&2240_x'&\cdot&\cdot&\cdot&\cdot&\cdot&\cdot&\cdot&\cdot&1&1&1&\cdot&\cdot\\
&3240_z'&\cdot&\cdot&\cdot&\cdot&\cdot&\cdot&\cdot&\cdot&\cdot&1&1&1&\cdot\\
&1400_z'&\cdot&\cdot&\cdot&\cdot&\cdot&\cdot&\cdot&\cdot&\cdot&\cdot&1&\cdot&\cdot\\
&1008_z'&\cdot&\cdot&\cdot&\cdot&\cdot&\cdot&\cdot&\cdot&\cdot&\cdot&\cdot&1&\cdot\\
&56_z'&\cdot&\cdot&\cdot&\cdot&\cdot&\cdot&\cdot&\cdot&\cdot&\cdot&\cdot&\cdot&1\\
\end{block}\\
\begin{block}{cc(ccccccccccccc)}
\mathtt{(4)}&8_z&1&\cdot&\cdot&-1&\cdot&1&-1&\cdot&\cdot&\cdot&\cdot&\cdot&1\\
\mathtt{(4)}&560_z&\cdot&1&-1&-1&1&1&-1&\cdot&1&-1&\cdot&1&\cdot\\
\mathtt{(5)}&1344_x&\cdot&\cdot&1&\cdot&-1&-1&1&\cdot&\cdot&1&-1&-1&\cdot\\
\mathtt{(5)}&1400_{zz}&\cdot&\cdot&\cdot&1&\cdot&-1&1&\cdot&-1&1&\cdot&-1&-1\\
\mathtt{(4)}&840_z&\cdot&\cdot&\cdot&\cdot&1&\cdot&-1&1&\cdot&-1&1&\cdot&\cdot\\
&4536_z&\cdot&\cdot&\cdot&\cdot&\cdot&1&-1&\cdot&\cdot&-1&1&1&1\\
&3200_x&\cdot&\cdot&\cdot&\cdot&\cdot&\cdot&1&-1&\cdot&1&-1&\cdot&-1\\
&4200_z'&\cdot&\cdot&\cdot&\cdot&\cdot&\cdot&\cdot&1&\cdot&-1&1&\cdot&\cdot\\
&2240_x'&\cdot&\cdot&\cdot&\cdot&\cdot&\cdot&\cdot&\cdot&1&-1&\cdot&1&\cdot\\
&3240_z'&\cdot&\cdot&\cdot&\cdot&\cdot&\cdot&\cdot&\cdot&\cdot&1&-1&-1&\cdot\\
&1400_z'&\cdot&\cdot&\cdot&\cdot&\cdot&\cdot&\cdot&\cdot&\cdot&\cdot&1&\cdot&\cdot\\
&1008_z'&\cdot&\cdot&\cdot&\cdot&\cdot&\cdot&\cdot&\cdot&\cdot&\cdot&\cdot&1&\cdot\\
&56_z'&\cdot&\cdot&\cdot&\cdot&\cdot&\cdot&\cdot&\cdot&\cdot&\cdot&\cdot&\cdot&1\\
\end{block}
\end{blockarray}
\]

All decomposition numbers not belonging to the Hecke algebra decomposition matrix follow immediately from Lemmas \ref{RR} and \ref{dim Hom} except for $[\m(560_z):\el(840_z)]$. Let $\alpha:=[\m(560_z):\el(840_z)]$ and write:
\begin{align*}
\el(560_z)&=\m(560_z)-\m(1344_x)-\m(1400_{zz})+(1-\alpha)\m(840_z)+\m(4536_z)+(\alpha-1)\m(3200_x)\\&-\alpha\m(4200_z')+\m(2240_x')+(\alpha-1)\m(3240_z')-\alpha\m(1400_z')+\m(1008_z')\\
\end{align*}
Since $\el(840_z')=\m(840_z')-\m(1344_x')+\m(560_z')$ and $\dim\Hom(\m(560_z'),\m(1344_x'))=1$, $\dim\Hom(\m(560_z),\m(840_z))=\dim\Hom(\m(560_z'),\m(840_z'))=0$,\\ $\dim\Hom(\m(840_z),\m(1400_{zz}))=0$, and $\dim\Hom(\m(840_z),\m(1344_x))=1$. Therefore $\alpha$ can equal $0$ or $1$. Suppose $\alpha=1$ and let $M$ be the supposed expression for $\el(560_z)$ in this case. Restricting $M$ to $H_{\frac{1}{4}}(E_7)$,
\begin{align*}
\Res M&=\m(7_a')+\m(56_a')-\m(210_a)-\m(105_b)-\m(280_b)-\m(15_a')+\m(405_a)+\m(378_a')\\
&-\m(105_c')-\m(315_a)-\m(216_a)+\m(210_b')+\m(189_b)+\m(70_a)+\m(21_a')+\m(120_a')+\m(7_a)\\
\end{align*}
Since $\m(7_a')$ and $\m(56_a')$ belong to different blocks and have minimal $h_c$-weights for those blocks over all the $\m(\sigma)$'s appearing in $\Res M$, $\el(7_a')$ and $\el(56_a')$ must appear in the composition series of $\Res M$. Subtracting them off we get that 
$$\Res M-\el(7_a')-\el(56_a')=-\m(189_a)+\m(336_a')-\m(315_a)+2\m(189_b)+\m(7_a)$$
and among the terms on the right-hand-side, $189_a$ has the lowest $h_c$-weight, so $-\m(189_a)$ can't come from the composition series of any of the other terms. Thus $\Res M$ can't be a module, $M$ can't either, $\alpha$ can't be $1$ and so $\alpha=0$.

And the block dual to that containing $\el(8_z)$:
\[
\begin{blockarray}{ccccccccccccccc}
\begin{block}{cc(ccccccccccccc)}
\mathtt{(4)}&1400_z&1&\cdot&\cdot&1&1&\cdot&\cdot&\cdot&\cdot&\cdot&\cdot&\cdot&\cdot\\
\mathtt{(4)}&1008_z&\cdot&1&\cdot&1&\cdot&1&1&\cdot&\cdot&\cdot&\cdot&\cdot&\cdot\\
\mathtt{(4)}&56_z&\cdot&\cdot&1&\cdot&\cdot&\cdot&1&\cdot&\cdot&1&\cdot&\cdot&\cdot\\
\mathtt{(5)}&3240_z&\cdot&\cdot&\cdot&1&1&1&\cdot&\cdot&1&\cdot&\cdot&\cdot&\cdot\\
&2240_x&\cdot&\cdot&\cdot&\cdot&1&\cdot&\cdot&\cdot&1&\cdot&\cdot&\cdot&1\\
&4200_z&\cdot&\cdot&\cdot&\cdot&\cdot&1&1&1&1&\cdot&\cdot&\cdot&\cdot\\
\mathtt{(5)}&3200_x'&\cdot&\cdot&\cdot&\cdot&\cdot&\cdot&1&1&\cdot&1&1&\cdot&\cdot\\
&4536_z'&\cdot&\cdot&\cdot&\cdot&\cdot&\cdot&\cdot&1&1&\cdot&1&1&\cdot\\
&1400_{zz}'&\cdot&\cdot&\cdot&\cdot&\cdot&\cdot&\cdot&\cdot&1&\cdot&\cdot&1&1\\
&840_z'&\cdot&\cdot&\cdot&\cdot&\cdot&\cdot&\cdot&\cdot&\cdot&1&1&\cdot&\cdot\\
&1344_x'&\cdot&\cdot&\cdot&\cdot&\cdot&\cdot&\cdot&\cdot&\cdot&\cdot&1&1&\cdot\\
&560_z'&\cdot&\cdot&\cdot&\cdot&\cdot&\cdot&\cdot&\cdot&\cdot&\cdot&\cdot&1&\cdot\\
&8_z'&\cdot&\cdot&\cdot&\cdot&\cdot&\cdot&\cdot&\cdot&\cdot&\cdot&\cdot&\cdot&1\\
\end{block}\\
\begin{block}{cc(ccccccccccccc)}
\mathtt{(4)}&1400_z&1&\cdot&\cdot&-1&\cdot&1&-1&\cdot&\cdot&1&\cdot&\cdot&\cdot\\
\mathtt{(4)}&1008_z&\cdot&1&\cdot&-1&1&\cdot&-1&1&-1&1&-1&1&\cdot\\
\mathtt{(4)}&56_z&\cdot&\cdot&1&\cdot&\cdot&\cdot&-1&1&-1&\cdot&\cdot&\cdot&1\\
\mathtt{(5)}&3240_z&\cdot&\cdot&\cdot&1&-1&-1&1&\cdot&1&-1&\cdot&-1&\cdot\\
&2240_x&\cdot&\cdot&\cdot&\cdot&1&\cdot&\cdot&\cdot&-1&\cdot&\cdot&1&\cdot\\
&4200_z&\cdot&\cdot&\cdot&\cdot&\cdot&1&-1&\cdot&-1&1&\cdot&1&1\\
\mathtt{(5)}&3200_x'&\cdot&\cdot&\cdot&\cdot&\cdot&\cdot&1&-1&1&-1&1&-1&-1\\
&4536_z'&\cdot&\cdot&\cdot&\cdot&\cdot&\cdot&\cdot&1&-1&\cdot&-1&1&1\\
&1400_{zz}'&\cdot&\cdot&\cdot&\cdot&\cdot&\cdot&\cdot&\cdot&1&\cdot&\cdot&-1&-1\\
&840_z'&\cdot&\cdot&\cdot&\cdot&\cdot&\cdot&\cdot&\cdot&\cdot&1&-1&1&\cdot\\
&1344_x'&\cdot&\cdot&\cdot&\cdot&\cdot&\cdot&\cdot&\cdot&\cdot&\cdot&1&-1&\cdot\\
&560_z'&\cdot&\cdot&\cdot&\cdot&\cdot&\cdot&\cdot&\cdot&\cdot&\cdot&\cdot&1&\cdot\\
&8_z'&\cdot&\cdot&\cdot&\cdot&\cdot&\cdot&\cdot&\cdot&\cdot&\cdot&\cdot&\cdot&1\\
\end{block}
\end{blockarray}
\]

All entries in columns not belonging to the Hecke algebra decomposition matrix follow immediately from Lemmas \ref{RR} and \ref{dim Hom} except for $[\m(1400_z):\el(3200_x')]$, \\$[\m(1008_z):\el(3200_x')]$, and $[\m(3240_z):\el(3200_x')]$.

Set $\alpha:=[\m(3240_z):\el(3200_x')]$, $\beta:=[\m(1008_z):\el(3200_x')]$, and \\$\gamma:=[\m(1400_z):\el(3200_x')]$. Write:
\begin{align*}
&\el(3240_z)=\m(3240_z)-\m(2240_x)-\m(4200_z)+(1-\alpha)\m(3200_x')+\alpha\m(4536_z')\\
&\qquad+(1-\alpha)\m(1400_{zz}')+(\alpha-1)\m(840_z')-\alpha\m(1344_x')+(\alpha-1)\m(560_z')+\alpha\m(8_z')\\
\\
&\el(1008_z)=\m(1008_z)-\m(3240_z)+\m(2240_x)-\beta\m(3200_x')+\beta\m(4536_z')-\beta\m(1400_{zz}')\\
&\qquad+\beta\m(840_z')-\beta\m(1344_x')+\beta\m(560_z')+(\beta-1)\m(8_z')\\
\\
&\el(1400_z)=\m(1400_z)-\m(3240_z)+\m(4200_z)-(\gamma+1)\m(3200_x')+\gamma\m(4536_z')-\gamma\m(1400_{zz}')\\
&\qquad+(\gamma+1)\m(840_z')-\gamma\m(1344_x')+\gamma\m(560_z')+\gamma\m(8_z')\\
\end{align*}

First we show that $\alpha=0$. From the decompositions of $\el(3200_x),$ $\el(2240_x'),$ and $\el(4200_z')$ and Lemma \ref{dim Hom}, it follows that $\dim\Hom(\m(3200_x'),\m(3240_z))=\dim\Hom(\m(3200_x'),\m(2240_x))=0,$ while $\dim\Hom(\m(3200_x'),\m(4200_z))=1$. Thus $\alpha$ is either $0$ or $1$. Suppose $\alpha=1$, so that $\el(3240_z)$ is equal to 
$$M:=\m(3240_z)-\m(2240_x)-\m(4200_z)+\m(4536_z')-\m(1344_x')+\m(8_z')$$
and calculate the restriction of $M$ to $H_{\frac{1}{4}}(E_7)$, and consider only those terms which belong to the block containing $\el(35_b)$:
$$\Res M|_{\B(35_b)}=\m(35_b)-\m(210_b)-\m(105_c)+\m(378_a)-\m(105_a)+\m(7_a)$$
Since $35_b$ has the smallest $h_c$-weight of any of the lowest weights of the Vermas in this expression, $\el(35_b)$ must belong to the composition series of $\Res M|_{\B(35_b)}$. Subtracting off $\el(35_b)=\m(35_b)-\m(210_b)+\m(280_b')-\m(105_a),$ we have:
$$\Res M|_{\B(35_b)}-\el(35_b)=-\m(105_c)+\m(378_a)-\m(280_b')+\m(7_a)$$
which cannot be a module in this block since the smallest $h_c$-weight term, $-\m(105_c)$, appears with a minus sign.
Therefore $\Res M|_{\B(35_b)}$ is not a module and so $M$ cannot be either, and it must be that $\alpha=0$.

Next, we claim that $\beta=1$. Consider
$$\Res^{E_8}_{E_7}1008_z=27_a\oplus21_a\oplus120_a\oplus210_a\oplus105_a'\oplus189_c'\oplus56_a'\oplus280_a'$$
Then observe that $\m(27_a)$ is the only term in $\el(27_a)$ with lowest weight one of the factors of $\Res^{E_8}_{E_7}1008_z$. Since $\m(1008_z)$ does not appear in the composition series of any simple in its block except for $\el(1008_z)$, it follows that $\el(1008_z)$ appears in the decomposition of $\Ind\el(27_a)$ into simples.  This implies that $\dim\Supp\el(1008_z)\leq4$, since $\dim\Supp\el(27_a)=3$ and induction from a maximal parabolic raises dimension of support of a module by $1$. On the other hand, the character of $\el(1008_z)$ will only have a pole of order $4$ if $\beta=1$, otherwise it will have a pole of order $5$. So $\beta$ must be $1$ and $\dim\Supp\el(1008_z)=4$.

The same kind of argument shows $\gamma=0$: the character of $\el(1400_z)$ has a pole of order $5$ if $\gamma\geq1$, but if $\gamma=0$, it has a pole of order $4$. But $\el(1400_z)$ appears in the composition series of $\Ind_{E_7}^{E_8}\el(21_b')$ into simples, and $\dim\Supp\el(21_b')=3$, so $\dim\Supp\el(1400_z)\leq4$ and $\gamma=0$.

\section{$F_4$ with equal parameters}

The relevant parameters are $c=1/d$ for $d=2,3,4,6,8,12$. We omit the case $d=2$. For $d=8,12$ the blocks are of defect $1$ (\cite{GP}, Table F.3) and the unique representation of less than full support in each block happens to be finite-dimensional -- the characters were originally found by Rouquier \cite{R}. All other nontrivial blocks have defect $2$ (\cite{GP}, Table F.3). 

\subsection{{\Large$\mbf{c=\frac{1}{6}}$}} There is a single block of non-zero defect (\cite{GJ}, Table 7.7).
\theorem The decomposition matrix of the principal block is as follows:
\[
\begin{blockarray}{ccccccccccccccc}
&&1_1&2_3&2_1&8_3&8_1&12&9_2&9_3&8_4&8_2&2_4&2_2&1_4\\
\begin{block}{cc(ccccccccccccc)}
\star&1_1&1&\cdot&\cdot&1&1&1&\cdot&\cdot&\cdot&\cdot&\cdot&\cdot&\cdot\\
\star&2_3&\cdot&1&\cdot&1&\cdot&\cdot&\cdot&1&\cdot&\cdot&\cdot&\cdot&\cdot\\
\star&2_1&\cdot&\cdot&1&\cdot&1&\cdot&1&\cdot&\cdot&\cdot&\cdot&\cdot&\cdot\\
\mathtt{(1)}&8_3&\cdot&\cdot&\cdot&1&\cdot&1&\cdot&1&\cdot&1&\cdot&\cdot&\cdot\\
\mathtt{(1)}&8_1&\cdot&\cdot&\cdot&\cdot&1&1&1&\cdot&1&\cdot&\cdot&\cdot&\cdot\\
&12&\cdot&\cdot&\cdot&\cdot&\cdot&1&\cdot&\cdot&1&1&\cdot&\cdot&1\\
&9_2&\cdot&\cdot&\cdot&\cdot&\cdot&\cdot&1&\cdot&1&\cdot&1&\cdot&\cdot\\
&9_3&\cdot&\cdot&\cdot&\cdot&\cdot&\cdot&\cdot&1&\cdot&1&\cdot&1&\cdot\\
&8_4&\cdot&\cdot&\cdot&\cdot&\cdot&\cdot&\cdot&\cdot&1&\cdot&1&\cdot&1\\
&8_2&\cdot&\cdot&\cdot&\cdot&\cdot&\cdot&\cdot&\cdot&\cdot&1&\cdot&1&1\\
&2_4&\cdot&\cdot&\cdot&\cdot&\cdot&\cdot&\cdot&\cdot&\cdot&\cdot&1&\cdot&\cdot\\
&2_2&\cdot&\cdot&\cdot&\cdot&\cdot&\cdot&\cdot&\cdot&\cdot&\cdot&\cdot&1&\cdot\\
&1_4&\cdot&\cdot&\cdot&\cdot&\cdot&\cdot&\cdot&\cdot&\cdot&\cdot&\cdot&\cdot&1\\
\end{block}
\\
\begin{block}{cc(ccccccccccccc)}
\star&1_1&1&\cdot&\cdot&-1&-1&1&1&1&-1&-1&\cdot&\cdot&1\\
\star&2_3&\cdot&1&\cdot&-1&\cdot&1&\cdot&\cdot&-1&\cdot&1&\cdot&\cdot\\
\star&2_1&\cdot&\cdot&1&\cdot&-1&1&\cdot&\cdot&\cdot&-1&\cdot&1&\cdot\\
\mathtt{(1)}&8_3&\cdot&\cdot&\cdot&1&\cdot&-1&\cdot&-1&1&1&-1&\cdot&-1\\
\mathtt{(1)}&8_1&\cdot&\cdot&\cdot&\cdot&1&-1&-1&\cdot&1&1&\cdot&-1&-1\\
&12&\cdot&\cdot&\cdot&\cdot&\cdot&1&\cdot&\cdot&-1&-1&1&1&1\\
&9_2&\cdot&\cdot&\cdot&\cdot&\cdot&\cdot&1&\cdot&-1&\cdot&\cdot&\cdot&1\\
&9_3&\cdot&\cdot&\cdot&\cdot&\cdot&\cdot&\cdot&1&\cdot&-1&\cdot&\cdot&1\\
&8_4&\cdot&\cdot&\cdot&\cdot&\cdot&\cdot&\cdot&\cdot&1&\cdot&-1&\cdot&-1\\
&8_2&\cdot&\cdot&\cdot&\cdot&\cdot&\cdot&\cdot&\cdot&\cdot&1&\cdot&-1&-1\\
&2_4&\cdot&\cdot&\cdot&\cdot&\cdot&\cdot&\cdot&\cdot&\cdot&\cdot&1&\cdot&\cdot\\
&2_2&\cdot&\cdot&\cdot&\cdot&\cdot&\cdot&\cdot&\cdot&\cdot&\cdot&\cdot&1&\cdot\\
&1_4&\cdot&\cdot&\cdot&\cdot&\cdot&\cdot&\cdot&\cdot&\cdot&\cdot&\cdot&\cdot&1\\
\end{block}
\end{blockarray}
\]

\begin{proof}

The $h_c$-weights on lowest weights in this block are given by $h_c(\tau)=2(1-\frac{\chi_\tau(s_1)+\chi_\tau(s_3)}{\dim\tau})$:\\

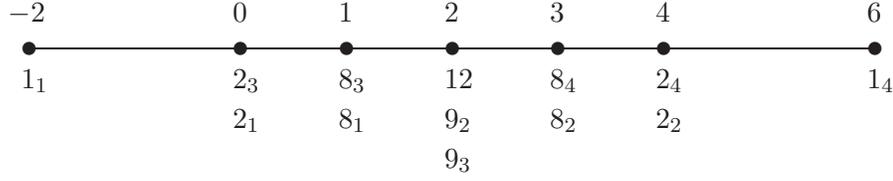
\begin{center}
\begin{picture}(320,70)
\put(0,50){\line(1,0){320}}
\put(0,50){\circle*{5}}
\put(80,50){\circle*{5}}
\put(120,50){\circle*{5}}
\put(160,50){\circle*{5}}
\put(200,50){\circle*{5}}
\put(240,50){\circle*{5}}
\put(320,50){\circle*{5}}
\put(-8,60){$-2$}
\put(77,60){$0$}
\put(117,60){$1$}
\put(157,60){$2$}
\put(197,60){$3$}
\put(237,60){$4$}
\put(317,60){$6$}
\put(-3,35){$1_1$}
\put(77,35){$2_3$}
\put(77,20){$2_1$}
\put(117,35){$8_3$}
\put(117,20){$8_1$}
\put(157,35){$12$}
\put(157,20){$9_2$}
\put(157,5){$9_3$}
\put(197,35){$8_4$}
\put(197,20){$8_2$}
\put(237,35){$2_4$}
\put(237,20){$2_2$}
\put(317,35){$1_4$}
\end{picture}
\end{center}

All but the first five columns of the decomposition matrix belong to that of the Hecke algebra. Of the first five columns, all entries are obtained from Lemmas \ref{RR} and \ref{dim Hom} except for $[\m(1_1):\el(8_3)]$ and $[\m(1_1):\el(8_1)]$. Since $\el(1_1)$ is finite-dimensional \cite{VV}, $\dim\el(1_1)[1]=\dim\el(1_1)[-1]$, and this implies $[\m(1_1):\el(8_3)]+[\m(1_1):\el(8_3)]=2$. Again because $\el(1_1)$ is finite-dimensional, writing $\el(1_1)=\sum c_\tau\m(\tau)$, we must have $\sum c_\tau\m(\tau')=\sum c_\tau\m(\tau)$; it follows that $[\m(1_1):\el(8_3)]=[\m(1_1):\el(8_1)]$.
\end{proof}

\subsection{{\Large$\mbf{c=\frac{1}{4}}$}} There are two defect $1$ blocks and the principal block, which has defect $2$ \cite{GJ} (Table 7.7). The two defect $1$ blocks with their $h_c$-weights are pictured below:

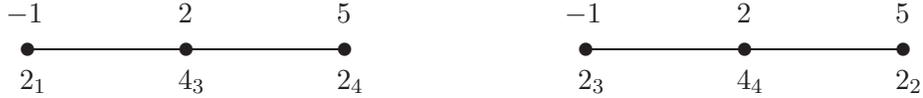
\begin{center}
\begin{picture}(120,40)
\put(0,20){\line(1,0){120}}
\put(0,20){\circle*{5}}
\put(60,20){\circle*{5}}
\put(120,20){\circle*{5}}
\put(-8,30){$-1$}
\put(57,30){$2$}
\put(117,30){$5$}
\put(-3,5){$2_1$}
\put(57,5){$4_3$}
\put(117,5){$2_4$}
\end{picture}
\qquad\qquad\qquad\qquad
\begin{picture}(120,40)
\put(0,20){\line(1,0){120}}
\put(0,20){\line(1,0){120}}
\put(0,20){\circle*{5}}
\put(60,20){\circle*{5}}
\put(120,20){\circle*{5}}
\put(-8,30){$-1$}
\put(57,30){$2$}
\put(117,30){$5$}
\put(-3,5){$2_3$}
\put(57,5){$4_4$}
\put(117,5){$2_2$}
\end{picture}
\end{center}

from which we calculate the characters of $\el(2_1)$ and $\el(2_3)$ and find that $$\dim\Supp\el(2_1)=\dim\Supp\el(2_3)=2$$

\theorem The decomposition matrix of the principal block is as follows:
\[
\begin{blockarray}{ccccccccccc}
&&1_1&4_2&9_1&12&6_1&4_1&9_4&4_5&1_4\\
\begin{block}{cc(ccccccccc)}
\star&1_1&1&\cdot&1&\cdot&1&\cdot&\cdot&\cdot&\cdot\\
\star&4_2&\cdot&1&1&1&\cdot&\cdot&\cdot&\cdot&\cdot\\
\mathtt{(2)}&9_1&\cdot&\cdot&1&1&1&1&1&\cdot&\cdot\\
&12&\cdot&\cdot&\cdot&1&\cdot&\cdot&1&1&\cdot\\
&6_1&\cdot&\cdot&\cdot&\cdot&1&\cdot&1&\cdot&1\\
\mathtt{(2)}&4_1&\cdot&\cdot&\cdot&\cdot&\cdot&1&1&\cdot&\cdot\\
&9_4&\cdot&\cdot&\cdot&\cdot&\cdot&\cdot&1&1&1\\
&4_5&\cdot&\cdot&\cdot&\cdot&\cdot&\cdot&\cdot&1&\cdot\\
&1_4&\cdot&\cdot&\cdot&\cdot&\cdot&\cdot&\cdot&\cdot&1\\
\end{block}
\\
\begin{block}{cc(ccccccccc)}
\star&1_1&1&\cdot&-1&1&\cdot&1&-1&\cdot&1\\
\star&4_2&\cdot&1&-1&\cdot&1&1&-1&1&\cdot\\
\mathtt{(2)}&9_1&\cdot&\cdot&1&-1&-1&-1&2&-1&-1\\
&12&\cdot&\cdot&\cdot&1&\cdot&\cdot&-1&\cdot&1\\
&6_1&\cdot&\cdot&\cdot&\cdot&1&\cdot&-1&1&\cdot\\
\mathtt{(2)}&4_1&\cdot&\cdot&\cdot&\cdot&\cdot&1&-1&1&1\\
&9_4&\cdot&\cdot&\cdot&\cdot&\cdot&\cdot&1&-1&-1\\
&4_5&\cdot&\cdot&\cdot&\cdot&\cdot&\cdot&\cdot&1&\cdot\\
&1_4&\cdot&\cdot&\cdot&\cdot&\cdot&\cdot&\cdot&\cdot&1\\
\end{block}
\end{blockarray}
\]

\begin{proof} The $h_c$-weights on lowest weights $\tau$ in this block are:

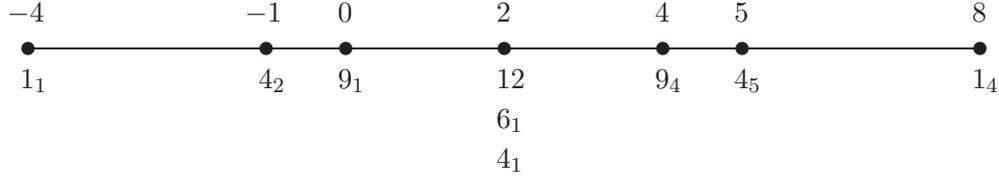
\begin{center}
\begin{picture}(360,70)
\put(0,50){\line(1,0){360}}
\put(0,50){\circle*{5}}
\put(90,50){\circle*{5}}
\put(120,50){\circle*{5}}
\put(180,50){\circle*{5}}
\put(240,50){\circle*{5}}
\put(270,50){\circle*{5}}
\put(360,50){\circle*{5}}
\put(-8,60){$-4$}
\put(82,60){$-1$}
\put(117,60){$0$}
\put(177,60){$2$}
\put(237,60){$4$}
\put(267,60){$5$}
\put(357,60){$8$}
\put(-3,35){$1_1$}
\put(87,35){$4_2$}
\put(117,35){$9_1$}
\put(177,35){$12$}
\put(177,20){$6_1$}
\put(177,5){$4_1$}
\put(237,35){$9_4$}
\put(267,35){$4_5$}
\put(357,35){$1_4$}
\end{picture}
\end{center}
The columns labeled by $\tau$ where $\el(\tau)$ has full support belong to the decomposition matrix of the Hecke algebra; of the four remaining columns, all entries follow from Lemmas \ref{RR} and \ref{dim Hom} except for $[\m(1_1):\el(4_1)]$ and $[\m(4_2):\el(4_1)]$. In order that dimensions of graded pieces of $\el(1_1)$ satisfy $\dim\el(1_1)[2]\geq\dim\el(1_1)[-2]$, $[\m(1_1):\el(4_1)]$ must be $0$. Likewise, in order to have $\dim\el(4_2)[2]\geq0$, $[\m(4_2):\el(4_1)]$ must be $0$.
\end{proof}

\subsection{{\Large$\mbf{c=\frac{1}{3}}$}} There are two blocks of defect $2$, one containing $\el(1_1)$ (the spherical representation) and the other containing $\el(4_2)$ ($4_2$ denotes the standard rep of $F_4$). Their decomposition matrices are identical, and if $\sigma_i,\;\tau_i\in\Irr F_4$ label the $i$'th rows in the decomposition matrices of the blocks of the trivial rep and standard rep respectively, then $\dim\tau_i=4\dim\sigma_i$. Moreover, the dimensions of supports of $\el(\tau_i)$ and $\el(\sigma_i)$ are the same.
 
\theorem The decomposition matrices of the blocks containing $\el(1_1)$ and $\el(4_2)$ are the same and given by the matrix below.

\[ 
\begin{blockarray}{cccccccccccc}
&&&4_2&8_3&8_1&16&4_3&4_4&8_4&8_2&4_5\\
&&&1_1&2_3&2_1&4_1&1_2&1_3&2_4&2_2&1_4\\
\begin{block}{ccc(ccccccccc)}
\star&4_2&1_1&1&1&1&1&\cdot&\cdot&\cdot&\cdot&\cdot\\
\mathtt{(2)}&8_3&2_3&\cdot&1&\cdot&1&\cdot&1&\cdot&1&\cdot\\
\mathtt{(2)}&8_1&2_1&\cdot&\cdot&1&1&1&\cdot&1&\cdot&\cdot\\
&16&4_1&\cdot&\cdot&\cdot&1&\cdot&\cdot&1&1&1\\
\mathtt{(2)}&4_3&1_2&\cdot&\cdot&\cdot&\cdot&1&\cdot&1&\cdot&\cdot\\
\mathtt{(2)}&4_4&1_3&\cdot&\cdot&\cdot&\cdot&\cdot&1&\cdot&1&\cdot\\
&8_4&2_4&\cdot&\cdot&\cdot&\cdot&\cdot&\cdot&1&\cdot&1\\
&8_2&2_2&\cdot&\cdot&\cdot&\cdot&\cdot&\cdot&\cdot&1&1\\
&4_5&1_4&\cdot&\cdot&\cdot&\cdot&\cdot&\cdot&\cdot&\cdot&1\\
\end{block}
\\
\begin{block}{ccc(ccccccccc)}
\star&4_2&1_1&1&-1&-1&1&1&1&-1&-1&1\\
\mathtt{(2)}&8_3&2_3&\cdot&1&\cdot&-1&\cdot&-1&1&1&-1\\
\mathtt{(2)}&8_1&2_1&\cdot&\cdot&1&-1&-1&\cdot&1&1&-1\\
&16&4_1&\cdot&\cdot&\cdot&1&\cdot&\cdot&-1&-1&1\\
\mathtt{(2)}&4_3&1_2&\cdot&\cdot&\cdot&\cdot&1&\cdot&-1&\cdot&1\\
\mathtt{(2)}&4_4&1_3&\cdot&\cdot&\cdot&\cdot&\cdot&1&\cdot&-1&1\\
&8_4&2_4&\cdot&\cdot&\cdot&\cdot&\cdot&\cdot&1&\cdot&-1\\
&8_2&2_2&\cdot&\cdot&\cdot&\cdot&\cdot&\cdot&\cdot&1&-1\\
&4_5&1_4&\cdot&\cdot&\cdot&\cdot&\cdot&\cdot&\cdot&\cdot&1\\
\end{block}
\end{blockarray}
\]

\begin{proof}
Of the five columns killed by KZ functor, Lemmas \ref{RR} and \ref{dim Hom} recover all their entries except $[\m(1_1):\el(1_2)]$ and $[\m(1_1):\el(1_3)]$ in the case of the principal block, correspondingly $[\m(4_2):\el(4_3)]$ and $[\m(4_2):\el(4_4)]$ in the other block. 

Consider the $h_c$-weights on the $\tau$ in each block:

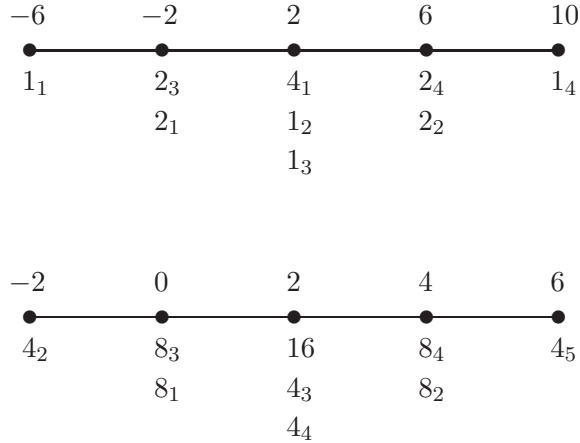
\begin{center}
\begin{picture}(200,80)
\put(0,60){\line(1,0){200}}
\put(0,60){\circle*{5}}
\put(50,60){\circle*{5}}
\put(100,60){\circle*{5}}
\put(150,60){\circle*{5}}
\put(200,60){\circle*{5}}
\put(-8,70){$-6$}
\put(42,70){$-2$}
\put(97,70){$2$}
\put(147,70){$6$}
\put(197,70){$10$}
\put(-3,45){$1_1$}
\put(47,45){$2_3$}
\put(47,30){$2_1$}
\put(97,45){$4_1$}
\put(97,30){$1_2$}
\put(97,15){$1_3$}
\put(147,45){$2_4$}
\put(147,30){$2_2$}
\put(197,45){$1_4$}
\end{picture}\\
\begin{picture}(200,100)
\put(0,60){\line(1,0){200}}
\put(0,60){\circle*{5}}
\put(50,60){\circle*{5}}
\put(100,60){\circle*{5}}
\put(150,60){\circle*{5}}
\put(200,60){\circle*{5}}
\put(-8,70){$-2$}
\put(47,70){$0$}
\put(97,70){$2$}
\put(147,70){$4$}
\put(197,70){$6$}
\put(-3,45){$4_2$}
\put(47,45){$8_3$}
\put(47,30){$8_1$}
\put(97,45){$16$}
\put(97,30){$4_3$}
\put(97,15){$4_4$}
\put(147,45){$8_4$}
\put(147,30){$8_2$}
\put(197,45){$4_5$}
\end{picture}
\end{center}

A dimension check on graded pieces of modules then forces the two missing decomposition numbers to be $0$ in both blocks:
\begin{align*}
&31=\dim\el(1_1)[-2]\leq\dim\el(1_1)[2]=29+1-[\m(1_1):\el(1_2)]+1-[\m(1_1):\el(1_3)]\\
\\
&4=\dim\el(4_2)[-2]\leq\dim\el(4_2)[2]=-4+4(1-[\m(4_2):\el(4_3)])+4(1-[\m(4_2):\el(4_4)])\\
\end{align*}
\end{proof}

\end{document}